\documentclass{article}

\usepackage{amssymb}
\usepackage{amsthm}
\usepackage{amscd}
\usepackage{amsmath}
\theoremstyle{plain}
\newtheorem{theorem}{Theorem}[section]
\newtheorem*{theorem*}{Theorem}

\newtheorem*{maintheorem*}{Main Theorem}
\newtheorem{proposition}[theorem]{Proposition}
\newtheorem{corollary}[theorem]{Corollary}
\newtheorem{lemma}[theorem]{Lemma}

\newtheorem*{conjecture*}{Conjecture}

\newtheorem{question}[theorem]{Question}
\theoremstyle{definition}
\newtheorem{definition}[theorem]{Definition}
\newtheorem*{definition*}{Definition}
\newtheorem{example}[theorem]{Example}
\newtheorem*{example*}{Example}

\newtheorem*{notation*}{Notation}
\newtheorem*{notation-conv*}{Notation and convention}
\newtheorem*{convention*}{Convention}
\theoremstyle{remark}
\newtheorem{remark}{Remark}
\newcommand{\Z}{{\mathbb Z}}
\newcommand{\C}{{\mathbb C}}
\newcommand{\Q}{{\mathbb Q}}
\newcommand{\R}{{\mathbb R}}
\newcommand{\proj}{{\mathbb P}}

\newcommand{\NN}{{\mathbb N}}

\begin{document}
\title{On the $0$-dimensional cusps of the K\"ahler moduli of a $K3$ surface} 
\author{Shouhei Ma\thanks{Graduate School of Mathematical Sciences, University of Tokyo  \endgraf
e-mail: sma@ms.u-tokyo.ac.jp } }
\date{}
\maketitle

\begin{abstract}
Let $S$ be a complex algebraic $K3$ surface. 
It is proved that 
the $0$-dimensional cusps of the K\"ahler moduli of $S$ 
are 
in one-to-one correspondence with 
the twisted Fourier-Mukai partners of $S$. 
As a result,    
a counting formula for the $0$-dimensional cusps of the K\"ahler moduli is obtained. 
Applications to 
rational maps between $K3$ surfaces 
are given.   
When the Picard number of $S$ is $1$, 
the bijective correspondence is calculated explicitly  
by using the Fricke modular curve.
\end{abstract}



\section{Introduction}\label{intro}

Let $S$ be an algebraic $K3$ surface over the field of complex numbers. 
In \cite{Ma1} 
we considered a certain orthogonal modular variety 
$\mathcal{K} (S) = \Gamma (S)^{+} \backslash \Omega _{\widetilde{NS}(S)}^{+}$ 
and studied the Fourier-Mukai (FM) partners of $S$ 
in connection with 
a compactification of 
$\mathcal{K} (S)$. 
It is a classical result of Baily-Borel \cite{B-B} 
that $\mathcal{K} (S)$ can be compactified 
to be a normal projective variety, 
by adjoining certain boundary components. 
The boundary components of  
$\mathcal{K} (S)$ 
are called {\it cusps\/} of 
$\mathcal{K} (S)$. 
One of the results of \cite{Ma1} 
is that 
the isomorphism classes of the FM partners of $S$ 
are in bijective correspondence with 
certain $0$-dimensional cusps of $\mathcal{K} (S)$, 
called {\it standard cusps\/}.       
However, 
a general $0$-dimensional cusp of $\mathcal{K} (S)$ 
is not necessarily standard. 
The purpose of this paper 
is to study the geometric counterparts 
of the {\it non-standard\/} $0$-dimensional cusps of $\mathcal{K} (S)$.

As we shall explain below, 
a non-standard cusp corresponds to 
a {\it twisted\/} Fourier-Mukai partner of $S$. 
The notion of twisted FM partners, 
which was introduced by C\u ald\u araru \cite{Ca1}, 
generalizes that of FM partners. 
By definition, 
a twisted FM partner of $S$ 
is a twisted $K3$ surface $(S', \alpha ')$ such that 
there is an exact equivalence 
$D^{b}(S', \alpha ') \simeq D^{b}(S)$ 
between the derived categories. 
As well as the FM partners, 
the twisted FM partners of $S$ can be regarded as 
certain geometric realizations of 
the category $D^{b}(S)$. 
They are useful for dealing with non-fine moduli spaces of sheaves on $S$ (see \cite{Ca2}, \cite{H-S1}).

Our main result is the following.

\begin{theorem}[Theorem \ref{most general}]\label{main1}
Let 
${\rm FM\/}^{d}(S) $ 
be the set of 
isomorphism classes of 
the twisted Fourier-Mukai partners $(S', \alpha ')$ of $S$ with ${\rm ord\/}(\alpha ')=d$.  
Let 
$\mathcal{C}^{d}(S)$ 
be the set of $0$-dimensional cusps of 
the Baily-Borel compactification of 
the modular variety  
$\mathcal{K} (S)$ 
of divisibility $d$. 
Then 
there exists a canonical bijection  
\begin{equation*}
{\rm FM\/}^{d}(S) 
 \simeq
\mathcal{C}^{d}(S). 
\end{equation*}
In particular,  we have 
\begin{equation*}
 \sum_{d} \# {\rm FM\/}^{d}(S) = 
\# \{ \text{ the 0-dimensional cusps of }  \mathcal{K} (S) \} . 
\end{equation*}
\end{theorem}

The existence of a map 
$\mathcal{C}^{d}(S) \to {\rm FM\/}^{d}(S)$ 
is observed in \cite{Ma1}. 
In the present paper 
we study the correspondence in more detail  
and establish its bijectivity.   
In fact,   
a more general result than Theorem \ref{main1}  will be proved.   
That is,  
for a {\it twisted\/} $K3$ surface $(S, \alpha )$  
we establish a relation between 
the $0$-dimensional cusps of the K\"ahler moduli of $(S, \alpha )$  
and the twisted FM partners of $(S, \alpha )$.    
If we put $d=1$ in Theorem \ref{main1}, 
we recover the correspondence between 
the untwisted FM partners 
and 
the standard cusps 
in \cite{Ma1}. 
Twisted FM partners appear naturally  
if we consider all $0$-dimensional cusps of $\mathcal{K} (S)$.

The modular variety $\mathcal{K} (S)$ is called K\"ahler moduli  
because it is an analogue of K\"ahler moduli of a Calabi-Yau $3$-fold.   
The analogy was investigated extensively by Dolgachev \cite{Do}.  
In an effort to formulate mirror symmetry for $K3$ surfaces, 
Dolgachev defined for a $K3$ surface $S$ (satisfying suitable conditions)  
a moduli space $\mathcal{M}^{\vee}$ of certain ``mirror $K3$ surfaces" of $S$,  
and observed that $\mathcal{M}^{\vee}$ is a modular variety 
uniformized by the Hermitian symmetric domain $\Omega _{\widetilde{NS}(S)}^{+}$ associated to $S$.   
Our K\"ahler moduli $\mathcal{K} (S)$  
is almost isomorphic to Dolgachev's mirror moduli space $\mathcal{M}^{\vee}$, 
in the sense that 
the arithmetic group $\Gamma (S)^{+}$ defining $\mathcal{K} (S)$ 
contains the arithmetic group defining $\mathcal{M}^{\vee}$ as a finite-index subgroup.   
When the Hodge structure of $S$ is generic, 
$\mathcal{K} (S)$ is indeed isomorphic to $\mathcal{M}^{\vee}$.   
We note that $\mathcal{K} (S)$ is naturally dominated by the complexified K\"ahler cone of $S$ (see \cite{Do}).    
On the other hand,  
Bridgeland gave an intrinsic construction of $\mathcal{K} (S)$ 
by proving that $\mathcal{K} (S)$ contains 
a natural quotient space of the space of stability conditions on the category $D^{b}(S)$ 
as a Zariski open set (see \cite{Br}).

An observation for the mirror picture of \cite{Do} can be drawn from Theorem \ref{main1}. 
In the formulation of \cite{Do},  
there is an ambiguity of the choice of mirror family,  
which depends on the choice of $0$-dimensional cusp of a moduli space $\mathcal{M}$ to which $S$ belongs.  
The moduli space $\mathcal{M}$ can be identified with 
the K\"ahler moduli of a generic member $S^{\vee}$ of a mirror moduli space of $S$.  
Theorem \ref{main1},  applied to $S^{\vee}$,  suggests that  
the ambiguity of the choice of mirror family of $S$ comes from the existence of twisted FM partners of $S^{\vee}$.

There is an application of Theorem \ref{main1}  
to the modular variety 
$\mathcal{K} (S)$.   
The twisted Fourier-Mukai number $\# {\rm FM\/}^{d}(S)$ is non-zero for only finitely many $d$,  
and each $\# {\rm FM\/}^{d}(S)$ is a finite number (see \cite{H-S1}).   
A formula for the number $\# {\rm FM\/}^{d}(S)$ is given in \cite{Ma2}. 
Roughly speaking, 
the number $\# {\rm FM\/}^{d}(S)$ is expressed 
as a sum of 
certain `masses' of the genera of some Lorentzian lattices.  
Combining those formulae with Theorem \ref{main1}, 
we obtain a counting formula 
for the $0$-dimensional cusps of 
the modular variety 
$\mathcal{K} (S)$ (Theorem \ref{number of cusps}). 
Via \cite{Do}, 
it allows us to 
count the $0$-dimensional cusps of the moduli spaces of 
certain lattice-polarized $K3$ surfaces.  
This formula for the cusps is a natural generalization of Scattone's formula \cite{Sc}.

When some description of the set of cusps is available, 
we obtain a classification of the twisted FM partners by Theorem \ref{main1}.  
We shall give such classification for certain elliptic $K3$ surfaces.

\begin{theorem}[Theorem \ref{twisted partners and relative Jacobians}]\label{main2}
Let $S$ be a $K3$ surface  
whose N\'eron-Severi lattice $NS(S)$ contains the hyperbolic plane $U$. 
Then 
for each twisted Fourier-Mukai partner $(S', \alpha ')$ of $S$ 
there exists an elliptic fibration $f : S \to  {\proj}^{1}$ such that
$(S', \alpha ')$ is isomorphic to 
the relative Jacobian of $f$. 
\end{theorem}

See Definition \ref{def of relative Jacobian} for the definition of relative Jacobian.  
In Kodaira's terminology, 
the elliptic surface underlying the relative Jacobian is 
the {\it basic elliptic surface\/} associated to the original fibration. 
The assumption that 
$NS(S)$ contains $U$ 
admits a geometric interpretation  
that $S$ has the structure of an elliptic surface with a section. 
For example, 
this assumption is satisfied 
if the Picard number $\rho (S)$ is larger than or equal to $13$.

Theorem \ref{main2} yields applications to rational maps between $K3$ surfaces.

\begin{corollary}[Proposition \ref{isogenous}]\label{main3} 
Let $S_{+}$ and $S_{-}$ be $K3$ surfaces 
with $\rho (S_{\pm}) \geq 13$. 
Let $T(S_{\pm})$ be the transcendental lattice of $S_{\pm}$. 
Then there exists a Hodge isometry $T(S_{+})\otimes {\Q} \simeq T(S_{-})\otimes {\Q}$ 
if and only if  
there exist a $K3$ surface $S_{0}$ and 
rational maps $S_{0} \dashrightarrow S_{+}$, $S_{0} \dashrightarrow S_{-}$  of square degrees. 
\end{corollary}

In other words,  
$S_{+}$ and $S_{-}$ are isogenous in the sense of Mukai \cite{Mu} 
if and only if they are dominated by a common $K3$ surface 
by rational maps of square degrees.

When $\rho (S) = 1$, 
we calculate the correspondence in Theorem \ref{main1} explicitly.  
Dolgachev \cite{Do} showed that 
the K\"ahler moduli $\mathcal{K} (S)$ for such $S$  is isomorphic to 
the {\it Fricke modular curve\/},  
which is the quotient of the congruence modular curve $\Gamma _{0}(n) \backslash {\mathbb H}$ 
by an involution. 
Here $2n$ is the degree of the $K3$ surface $S$. 
It is easy to describe the cusps of the Fricke curve,  
because the theory of elliptic modular curves is available.    
On the other hand,  
a set of representatives of ${\rm FM\/}^{d}(S)$ 
is given in \cite{Ma2} by certain moduli spaces of sheaves on $S$, 
twisted by natural obstruction classes.  
Then 
we have an explicit correspondence between 
the cusps of the Fricke curve 
and certain moduli spaces of sheaves on $S$. 
An advantage of 
considering the Fricke curve is that 
we have a complete description of its cusps.    
Thus, given two arbitrary primitive isotropic Mukai vectors, 
one can decide immediately  
whether the associated moduli spaces   
are isomorphic or not.

As we noted above,  
there is a derived-categorical construction of the K\"ahler moduli $\mathcal{K} (S)$ by using stability conditions due to Bridgeland. 
Then the following question arises naturally from Theorem \ref{main1}:

\begin{question}
Can one perform the Baily-Borel compactification of $\mathcal{K} (S)$ 
by studying degenerations of stability conditions in various large volume limits?  
\end{question}

The rest of the paper is structured as follows. 
In Sect.\ref{sec2.1} 
we recall some lattice theory. 
In Sect.\ref{sec2.2} 
we study lattice-theoretic properties of 
the twisted Mukai lattice of a twisted $K3$ surface. 
In Sect.\ref{sec3} 
we prove Theorem \ref{main1}.
In Sect.\ref{sec4} 
we derive Theorem \ref{main2} and Corollary \ref{main3}.
In Sect.\ref{sec5.1} 
we exhibit the isomorphism between 
the Fricke modular curve and the K\"ahler moduli of a $K3$ surface of Picard number $1$.
In Sect.\ref{sec5.2} 
we calculate the correspondence between 
the cusps of the Fricke modular curve and certain moduli spaces of sheaves.

\begin{notation*}
All varieties are assumed to be algebraic varieties over the field of complex numbers.  
In particular,  
a {\it $K3$ surface\/}   
means a nonsingular complex algebraic $K3$ surface.  
For a $K3$ surface $S$, 
we denote by $NS(S)$ (resp. $T(S)$) 
the N\'eron-Severi (resp. transcendental) lattice 
of $S$. 
The {\it Picard number\/} $\rho (S)$ 
is the rank of $NS(S)$. 
Let 
\begin{eqnarray*}
& & \widetilde{H}(S, {\Z}) 
         := H^{0}(S, {\Z}) + H^{2}(S, {\Z}) + H^{4}(S, {\Z}), \\
& & \widetilde{NS}(S) 
         := H^{0}(S, {\Z}) + NS(S) + H^{4}(S, {\Z}), 
\end{eqnarray*}
which are endowed with the Mukai pairing.  
The {\it hyperbolic plane\/} $U$ 
is the lattice ${\Z}e + {\Z}f$, 
$(e, e) = (f, f) =0, (e, f)=1$. 
We identify the lattice 
$H^{0}(S, {\Z}) + H^{4}(S, {\Z})$ 
endowed with the Mukai pairing 
with the hyperbolic plane 
$U$  
by the identifications  
$(1, 0, 0) = e$, 
$(0, 0, -1) = f$. 
Write 
\begin{equation*}
\Lambda _{K3} := E_{8}^{2} \oplus U^{3},  \: \: \: \: 
\widetilde{\Lambda }_{K3} := E_{8}^{2} \oplus U^{4}. 
\end{equation*}
\end{notation*}

\noindent{\textbf{Acknowledgements.}}
The author is deeply indebted to 
         Professor Ken-Ichi Yoshikawa 
                   for his encouragement and discussion on the whole part of the paper; 
         Professor Shinobu Hosono 
                   for discussion that motivated the work; 
         Doctor Kenji Hashimoto 
                   for discussion about quaternion orders.  
This work was partially supported by Grant-in-Aid for JSPS fellows.

\section{Twisted Mukai lattices}\label{sec2}

\subsection{Preliminaries from lattice theory}\label{sec2.1}

By an {\it even lattice\/}, 
we mean a free ${\Z}$-module $L$ of finite rank 
endowed with a non-degenerate symmetric bilinear form 
$L \times L \to {\Z}$ 
satisfying $(l, l) \in 2{\Z}$ 
for all $l \in L$. 
For a field $K$, 
the quadratic space $L\otimes K$ is denoted by $L_{K}$. 
For a vector $l \in L$ 
we define the {\it divisibility\/} 
${\rm div\/}(l)$ of $l$ 
as the positive generator of the ideal 
$(l, L) \subset {\Z}$. 
The set of the primitive isotropic vectors $l \in L$ 
with ${\rm div\/}(l) = d$ 
is denoted by $I^{d}(L)$. 
To an even lattice $L$ 
we can associate a finite Abelian group 
$D_{L} := L^{\vee} / L$ 
and a quadratic form 
$q_{L} : D_{L} \to {\Q}/ 2{\Z}$ 
defined by 
$q_{L} (x) \equiv (x, x) \: {\rm mod\/} \: 2{\Z}$, 
$x \in D_{L}$. 
Then 
$(D_{L}, q_{L})$ is called 
the {\it discriminant form\/} of $L$. 
We have a natural homomorphism 
$r_{L} : O(L) \to O(D_{L})$, 
whose kernel is denoted by $O(L)_{0}$. 
The following facts due to Nikulin \cite{Ni} are well-known. 
For later use, we indicate a proof.

\begin{proposition}[\cite{Ni}]\label{Nikulin}
Let $L$ be an even unimodular lattice 
and let $M$ be a primitive sublattice of $L$ 
with the orthogonal complement $M^{\perp}$.

$(1)$ 
There exists a natural isometry 
$\lambda _{L} : (D_{M}, q_{M}) \simeq (D_{M^{\perp}}, -q_{M^{\perp}})$.

$(2)$  
For two isometries 
$\gamma _{M} \in O(M)$ 
and 
$\gamma _{M^{\perp}} \in O(M^{\perp})$, 
$\gamma _{M} \oplus \gamma _{M^{\perp}}$ 
extends to the isometry of $L$ 
if and only if 
$r_{M}(\gamma _{M}) = 
\lambda _{L}^{-1} \circ r_{M^{\perp}}(\gamma _{M^{\perp}}) \circ \lambda _{L}$. 
\end{proposition}

\begin{proof}
$(1)$  
The unimodularity of $L$ and 
the primitivity of $M$  
assure the surjectivity of the projection 
$L \to M^{\vee}$. 
Thus 
for each $x \in M^{\vee} $ 
there is $y \in (M^{\perp})^{\vee}$ 
such that 
$x + y \in L$. 
The assignment $x \mapsto y$ 
induces the group isomorphism 
$\lambda _{L} : D_{M} \simeq D_{M^{\perp}}$. 
As  
$q_{M}(x) + q_{M^{\perp}}(y) \equiv (x+y, x+y) \in 2{\Z}$, 
we have 
$\lambda _{L} : (D_{M}, q_{M}) \simeq (D_{M^{\perp}}, -q_{M^{\perp}})$.

$(2)$ 
The isometry 
$\gamma _{M} \oplus \gamma _{M^{\perp}}$ 
of 
$M \oplus M^{\perp}$ 
preserves the overlattice $L$ 
if and only if 
$r_{M}(\gamma _{M}) \oplus r_{M^{\perp}}(\gamma _{M^{\perp}}) 
\in 
O(D_{M} \oplus D_{M^{\perp}})$ 
preserves the graph of 
$\lambda _{L}$. 
\end{proof}

Let $L$ be an even lattice of 
${\rm sign\/}(L) = (2, {\rm rk\/}(L)-2)$ 
and let 
$\Gamma \subset O(L)$ be a subgroup containing 
$\{ \pm {\rm id\/} \}$. 
Denote by  
$\Omega _{L}$ 
the set of the oriented positive-definite two-planes in $L_{{\R}}$,  
which carries a complex structure 
via the isomorphism 
\begin{equation*}\label{period domain and Grassmannian}
\Omega _{L} \: \simeq \: 
   \Bigl\{ 
           \Bigl. \: {\C}\omega \in {\proj}(L_{{\C}}) \: \Bigr| \: 
           (\omega , \omega )=0,  \: 
           (\omega , \bar{\omega})>0 \: 
   \Bigr\} . 
\end{equation*}
Then $\Omega _{L}$ has two connected components. 
A choice of a component of $\Omega _{L}$, 
say $\Omega _{L}^{+}$, 
is equivalent to 
a choice of an orientation for a positive-definite two-plane, 
and is called an {\it orientation\/} of $L$. 
Let 
$\Gamma ^{+} \subset \Gamma $ 
be the subgroup of 
the orientation-preserving isometries in $\Gamma $. 
The quotient space 
$\Gamma ^{+} \backslash \Omega _{L}^{+}$ 
admits the Baily-Borel compactification 
$\overline{\Gamma ^{+} \backslash \Omega _{L}^{+}}$, 
which turns out to be a normal projective variety (\cite{B-B}. See also \cite{Sc}). 
The set of 
the $0$-dimensional cusps of 
$\Gamma ^{+} \backslash \Omega _{L}^{+}$ 
is canonically identified with the set 
\begin{equation*}
\mathop{\bigcup}_{d} \: \Gamma ^{+} \backslash I^{d}(L). 
\end{equation*}

\subsection{Twisted Mukai lattices}\label{sec2.2}

The twisted Mukai lattice of a twisted $K3$ surface 
was defined by Huybrechts in \cite{Hu} 
and is studied in \cite{H-S1}. 
Here  
we develop lattice-theoretical properties of twisted Mukai lattices. 
The results of this section will be used in Section \ref{sec3}.

Let $S$ be a $K3$ surface. 
The {\it Brauer group\/} ${\rm Br\/}(S)$ of $S$ 
is the group of the torsion elements of 
$H^{2}(\mathcal{O}_{S}^{\times})$. 
Via the exponential sequence, 
we have 
\begin{equation}\label{eqn:Brauer group lift}
{\rm Br\/}(S) 
\simeq 
H^{2}(S, {\Q}) / 
( NS(S)_{{\Q}} + H^{2}(S, {\Z})) . 
\end{equation} 
For example, 
let $\rho (S) =20$. 
Then ${\rm Br\/}(S)$ is the group of the finite-order points of 
the elliptic curve 
$H^{2}(\mathcal{O}_{S}) / T(S)^{\vee}$.

A class $B\in H^{2}(S, {\Q})$ 
is called a {\it (rational) B-field lift\/} 
of 
a Brauer element $\alpha \in {\rm Br\/}(S)$ 
if $B$ maps to $\alpha $ in the isomorphism $(\ref{eqn:Brauer group lift})$. 
For an element $\alpha \in {\rm Br\/}(S)$ 
we can find a B-field lift of $\alpha $ 
from 
$\frac{1}{d}H^{2}(S, {\Z})$, 
where 
$d={\rm ord\/}(\alpha )$. 
By considering the intersection pairings 
of the B-field lifts 
with $T(S)$, 
we also have 
\begin{equation}\label{eqn:Brauer group hom}
{\rm Br\/}(S) 
\simeq 
{\rm Hom\/}(T(S), {\Q}/{\Z}). 
\end{equation}
Via $(\ref{eqn:Brauer group hom})$, 
we identify an element 
$\alpha \in {\rm Br\/}(S)$ 
with a surjective homomorphism 
$\alpha : T(S) \twoheadrightarrow {\Z}/{\rm ord\/}(\alpha ){\Z}$. 
Then 
${\rm Ker\/}(\alpha ) \subset T(S)$ 
is denoted by 
$T(S, \alpha )$.


A {\it twisted\/} $K3$ {\it surface\/} 
is a pair $(S, \alpha )$, 
where $S$ is a $K3$ surface 
and $\alpha \in {\rm Br\/}(S)$. 
Let 
$\omega _{S} \in H^{2}(S, {\C})$ 
be a period of $S$ 
and let $B \in H^{2}(S, {\Q})$ 
be a B-field lift of $\alpha $. 
By definition, the 
{\it twisted Mukai lattice\/} 
$\widetilde{H}(S, B, {\Z})$ of 
$(S, \alpha )$ and $B$ 
is the lattice $\widetilde{H}(S, {\Z})$ 
equipped with the twisted period 
${\rm e\/}^{B}(\omega _{S}) = 
(1, B, \frac{1}{2}(B, B)) \wedge (0, \omega _{S}, 0)$. 
Set  
\begin{eqnarray*}\label{twisted Neron-Severi and transcendental}
& & \widetilde{NS}(S, B) := 
                  {\rm e\/}^{B}(\omega _{S})^{\perp} \cap \widetilde{H}(S, B, {\Z}),  \\
& & T(S, B) := 
            \widetilde{NS}(S, B)^{\perp} \cap \widetilde{H}(S, B, {\Z}).
\end{eqnarray*} 
We have a Hodge isometry 
$
{\rm e\/}^{B} : T(S, \alpha ) \simeq T(S, B), 
$
as each class $l \in T(S)$ 
is of pure degree $2$ (\cite{H-S1}). 
Since 
$\widetilde{NS}(S, B)_{{\Q}} = {\rm e\/}^{B}(\widetilde{NS}(S)_{{\Q}})$, 
the orientation of 
$\widetilde{NS}(S)$ 
induces that of 
$\widetilde{NS}(S, B)$. 
That is, 
${\rm e\/}^{B}({\R}(1, 0, -1) \oplus {\R}(0, l, 0))$ 
is of positive orientation 
for an ample class  
$l \in NS(S)$.

For another B-field lift 
$B' \in H^{2}(S, {\Q})$ of $\alpha $, 
we can write 
$B' = B + B_{1} + B_{2}$ 
with 
$B_{1} \in H^{2}(S, {\Z})$ 
and 
$B_{2} \in NS(S)_{{\Q}}$. 
Then 
we have an orientation-preserving Hodge isometry 
\begin{equation*}
{\rm e\/}^{B_{1}}  : 
\widetilde{H}(S, B, {\Z}) 
\stackrel{\simeq}{\longrightarrow} 
\widetilde{H}(S, B+B_{1}, {\Z})=\widetilde{H}(S, B', {\Z})  
\end{equation*}
defined by the wedge product with the cohomology class 
$(1, B_{1}, \frac{1}{2}(B_{1}, B_{1}))$,   
which 
fixes the vector $(0, 0, 1)$.

In the remainder of this section, 
we fix a twisted $K3$ surface $(S, \alpha )$ 
and a B-field lift $B$ of $\alpha $.  
The basic ideas of the following Lemma \ref{twist and disc form} and Lemma \ref{div(0, 0, 1)} 
are present in \cite{Mu}. 

\begin{lemma}[\cite{Ma2}] \label{twist and disc form}
Let 
$\lambda : (D_{\widetilde{NS}(S, B)}, q) \simeq (D_{T(S, B)}, -q)$ 
be the natural isometry. 
Then we have a Hodge isometry 
\begin{equation}\label{eqn:extension of trans lattice}
{\rm e\/}^{B} : 
T(S) 
\stackrel{\simeq}{\to} 
\Bigl\langle T(S,B), \: \lambda ((0, 0, -\frac{1}{d})) \Bigr\rangle \subset T(S, B)^{\vee} , 
\end{equation}
where 
$d={\rm ord\/}(\alpha )$. 
The twisting 
$\alpha : T(S) \to {\Z}/d{\Z}$ 
is given by the homomorphism 
\begin{equation}\label{eqn:twist and disc form} 
T(S) 
\stackrel{{\rm e\/}^{B}}{\to} 
\Bigl\langle T(S,B), \: \lambda ((0, 0, -\frac{1}{d})) \Bigr\rangle / T(S, B) 
\simeq 
\Bigl\langle \lambda ((0, 0, -\frac{1}{d})) \Bigr\rangle 
\simeq 
{\Z}/d{\Z} .
\end{equation}
\end{lemma}

\begin{proof}
For a transcendental class 
$l\in T(S)$ 
with 
$\alpha (l) =\bar{1} \in {\Z}/d{\Z}$, 
we have 
${\rm e\/}^{B}(l) = (0, l, \frac{1}{d}+k)$ 
with 
$k \in {\Z}$. 
Since 
${\rm e\/}^{B}(l)+(0, 0, -\frac{1}{d}) \in \widetilde{H}(S, B, {\Z})$, 
we have 
$\lambda ((0, 0, -\frac{1}{d})) \equiv {\rm e\/}^{B}(l) \in D_{T(S, B)}$ 
by the definition of $\lambda $ (Proposition \ref{Nikulin}).  
Thus 
we obtain 
${\rm e\/}^{B}(T(S)) = 
\langle T(S,B), \: \lambda ((0, 0, -\frac{1}{d})) \rangle $. 
The image of $l$ 
by $(\ref{eqn:twist and disc form})$ 
is $\bar{1}$ 
and the image of $T(S, \alpha )$ 
by $(\ref{eqn:twist and disc form})$ 
is $\bar{0}$, 
so the second claim is proved. 
\end{proof}

\begin{lemma}\label{div(0, 0, 1)}
The divisibility of 
the primitive isotropic vector 
$(0, 0, 1) \in \widetilde{NS}(S, B)$
is equal to $d={\rm ord\/}(\alpha )$. 
\end{lemma}

\begin{proof}
Let 
$d' := {\rm div\/}((0, 0, 1))$. 
In the proof of Lemma \ref{twist and disc form} 
we saw that 
$(0, 0, \frac{1}{d}) \in \widetilde{NS}(S, B)^{\vee}$, 
which means that 
$d|d'$. 
Assume that $d<d'$. 
There exists a vector $l' \in T(S,B)^{\vee}$ 
such that 
$l' + (0, 0, \frac{1}{d'}) \in \widetilde{H}(S, B, {\Z})$. 
Since 
the primitive hull of 
$T(S, B) \oplus {\Z}(0, 0, 1)$ 
in 
$\widetilde{H}(S, B, {\Z})$ 
is 
$T(S) \oplus {\Z}(0, 0, 1)$, 
we can write 
$l' + (0, 0, \frac{1}{d'}) = l + (0, 0, k)$ 
for some 
$l \in T(S)$ and $k \in {\Z}$. 
Writing 
$l=l'' + (0, 0, \frac{k''}{d})$ 
with
$l'' := {\rm e\/}^{B}(l) \in T(S, B)^{\vee}$ and $k'' \in {\Z}$, 
we obtain the equality 
\begin{equation}\label{eqn:contradiction}
d'(l'-l'') = (0, 0, -1+\frac{d'}{d}k''+d'k).
\end{equation}
Then 
the right hand side of $(\ref{eqn:contradiction})$ 
is not $0$, 
which contradicts to the fact that 
$T(S, B)_{{\Q}} \cap {\Q}(0, 0, 1) = \{ 0 \}$. 
\end{proof}

Next we compare the two lattices 
$\widetilde{NS}(S)$ and $\widetilde{NS}(S, B)$. 
Consider the following isometry: 
\begin{equation}\label{kappa} 
\kappa : 
\widetilde{NS}(S) 
\stackrel{\simeq}{\to} 
dH^{0}(S, {\Z}) + NS(S) + \frac{1}{d}H^{4}(S, {\Z}), \: \: \: \: 
(a, l, b) \mapsto (da, l, \frac{1}{d}b), 
\end{equation}
where 
$d = {\rm ord\/}(\alpha )$.

\begin{proposition}\label{comparing}
When 
$B$ is chosen from $\frac{1}{d}H^{2}(S, {Z})$, 
$d = {\rm ord\/} (\alpha )$,  
we have the following orientation-preserving isometry: 
\begin{equation}\label{eqn:compare two NS}
{\rm e\/}^{B} \circ \kappa : 
\widetilde{NS}(S) \simeq 
\Bigl\langle \widetilde{NS}(S, B), \: (0, 0, -\frac{1}{d}) \Bigr\rangle . 
\end{equation}
\end{proposition}

\begin{proof}
Write 
$\widetilde{M} := \langle \widetilde{NS}(S, B), \: (0, 0, -\frac{1}{d}) \rangle $. 
From the equalities  
${\rm e\/}^{B}((0, 0, \frac{1}{d})) = (0, 0, \frac{1}{d})$, 
${\rm e\/}^{B}((d, 0, 0))           = (d, dB, \frac{1}{2}(dB, B))$, 
${\rm e\/}^{B}((0, l, 0))           = (0, l, (l, B))$, 
we obtain the inclusion
${\rm e\/}^{B} \circ \kappa (\widetilde{NS}(S)) \subset \widetilde{M}$. 
Since 
\begin{equation*}
{\rm det\/}\widetilde{M} 
= d^{-2}\cdot {\rm det\/}\widetilde{NS}(S, B) 
= d^{-2}\cdot {\rm det\/}T(S, \alpha ) 
= {\rm det\/}\widetilde{NS}(S), 
\end{equation*}
we have 
${\rm e\/}^{B} \circ \kappa (\widetilde{NS}(S)) = \widetilde{M}$.
\end{proof}

After choosing 
a B-field lift 
$B \in \frac{1}{d}H^{2}(S, {\Z})$ 
of $\alpha $, 
we write 
\begin{equation*}
\widetilde{M} := \langle \widetilde{NS}(S, B), \: (0, 0, -\frac{1}{d}) \rangle  , \: \: \: \: 
\widetilde{T} := \langle T(S, B), \: \lambda ((0, 0, -\frac{1}{d})) \rangle ,  
\end{equation*}
where 
$\lambda : D_{\widetilde{NS}(S, B)} \simeq D_{T(S, B)}$ 
is the natural isomorphism. 
Both $\widetilde{M}$ and $\widetilde{T}$ 
are even lattices. 
From $\lambda $ 
we obtain the isometry 
$\bar{\lambda } : (D_{\widetilde{M}}, q) \simeq (D_{\widetilde{T}}, -q)$. 
Denote by 
$\lambda _{0} : D_{\widetilde{NS}(S)} \simeq D_{T(S)}$ 
the natural isomorphism.

\begin{proposition}\label{commute}
Let $B \in \frac{1}{d}H^{2}(S, {\Z})$, $d= {\rm ord\/}(\alpha )$.  
The following diagram commutes. 
$$\CD
D_{\widetilde{NS}(S)} @>{\rm e\/}^{B} \circ \kappa >> D_{\widetilde{M}} \\
@V\lambda _{0}VV   @VV\bar{\lambda} V \\
D_{T(S)} @>>{\rm e\/}^{B}> D_{\widetilde{T}}.
\endCD $$
\end{proposition}

\begin{proof}
Every element of 
$D_{\widetilde{NS}(S)}$ 
can be represented by a vector in 
$NS(S)^{\vee}$. 
Thus 
we take a vector $x \in NS(S)^{\vee}$ 
and prove the commutativity for the element 
$[x] \in D_{NS(S)} \simeq D_{\widetilde{NS}(S)}$. 
Choose a vector 
$y \in T(S)^{\vee}$ 
representing 
$\lambda _{0}([x]) \in D_{T(S)}$. 
By the definition of $\lambda _{0}$, 
we have 
$x+y \in H^{2}(S, {\Z})$. 
Since 
${\rm e\/}^{B}(x+y) \in H^{2}(S, {\Z}) \oplus {\Z}(0, 0, \frac{1}{d})$, 
there is an integer 
$k \in {\Z}$ 
such that 
\begin{equation*}
{\rm e\/}^{B}(x) + (0, 0, \frac{k}{d}) + {\rm e\/}^{B}(y) 
\in H^{2}(S, {\Z}) 
\subset \widetilde{H}(S, B, {\Z}). 
\end{equation*}
By Lemma \ref{twist and disc form} 
and Proposition \ref{comparing}, 
we have 
${\rm e\/}^{B}(y) \in \widetilde{T}^{\vee}$ 
and 
${\rm e\/}^{B}(x) \in \widetilde{M}^{\vee}$. 
As 
$(0, 0, \frac{k}{d})$ is an integral vector in $\widetilde{M}$, 
we obtain 
${\rm e\/}^{B}(x) \equiv {\rm e\/}^{B}(x) + (0, 0, \frac{k}{d}) \in D_{\widetilde{M}}$. 
It follows that 
\[
 \bar{\lambda}({\rm e\/}^{B} \circ \kappa ([x]))
= \bar{\lambda}([{\rm e\/}^{B}(x)]) 
= \bar{\lambda}([{\rm e\/}^{B}(x) + (0, 0, \frac{k}{d})]) 
= [{\rm e\/}^{B}(y)] 
= {\rm e\/}^{B} \circ \lambda _{0} ([x]). 
\]
\end{proof}

\section{Twisted Fourier-Mukai partners and $0$-dimensional cusps}\label{sec3}

Let 
$(S, \alpha )$ 
be a twisted $K3$ surface 
with a B-field lift 
$B\in H^{2}(S, {\Q})$. 
A twisted $K3$ surface $(S', \alpha ')$ 
is called 
a {\it twisted Fourier-Mukai (FM) partner\/} of $(S, \alpha )$  
if there is an exact equivalence 
$D^{b}(S, \alpha ) \simeq D^{b}(S', \alpha ')$. 
Let 
${\rm FM\/}^{d}(S, \alpha )$ 
be the set of isomorphism classes of 
the twisted FM partners $(S', \alpha ')$ 
of $(S, \alpha )$ 
with ${\rm ord\/}(\alpha ')=d$.

\begin{definition}
We define the group 
$\Gamma (S, B) \subset O(\widetilde{NS}(S, B))$ by 
\begin{equation*}
\Gamma (S, B) := 
r_{\widetilde{NS}(S, B)}^{-1} (\lambda \circ r_{T(S, B)}(O_{Hodge}(T(S, B)))), 
\end{equation*}
where 
$\lambda : O(D_{T(S, B)}) \simeq O(D_{\widetilde{NS}(S, B)})$ 
is the isomorphism 
induced from the isometry 
$(D_{\widetilde{NS}(S, B)}, q) \simeq (D_{T(S, B)}, -q)$, 
and 
$O_{Hodge}(T(S, B))$
is the group of the Hodge isometries of $T(S, B)$. 
\end{definition}

There are obvious inclusions 
\begin{equation*}\label{eqn:inclusions} 
\{ \pm {\rm id\/} \} \times O(\widetilde{NS}(S, B))_{0} 
\subset 
\Gamma (S, B)
\subset 
O(\widetilde{NS}(S, B)). 
\end{equation*}
Each isometry 
$\gamma \in \Gamma (S, B)$ 
can be extended to a Hodge isometry 
$\widetilde{H}(S, B, {\Z}) \simeq \widetilde{H}(S, B, {\Z})$. 
Recall that 
$\Gamma (S, B)^{+}$ 
is the subgroup of $\Gamma (S, B)$ 
consisting of the orientation-preserving isometries in $\Gamma (S, B)$. 
Then 
we can form the modular variety 
\begin{equation}\label{eqn:twisted Kahler moduli}
\mathcal{K} (S, \alpha ) := 
          \Gamma (S, B)^{+} \backslash \Omega _{\widetilde{NS}(S, B)}^{+},  
\end{equation}  
the isomorphism class of which 
does not depend on the choices of the lift $B$ of $\alpha $. 
The set of $0$-dimensional cusps of the Baily-Borel compactification of $\mathcal{K} (S, \alpha )$ 
is identified with the set 
$\mathop{\bigcup}_{d} \; \mathcal{C}^{d}(S, \alpha )$, 
where 
\begin{equation}\label{eqn:0-dim cusps of twisted Kahler moduli}
\mathcal{C}^{d}(S, \alpha ) \: := \: 
                                \Gamma (S, B)^{+} \backslash I^{d}(\widetilde{NS}(S, B)). 
\end{equation}
When 
$\alpha =1$, 
we write 
$\mathcal{K} (S) := \mathcal{K} (S, 1)$
and 
$\mathcal{C}^{d} (S) := \mathcal{C}^{d} (S, 1)$. 
Then we have 
\begin{equation}\label{eqn:Kahler moduli}
\mathcal{K}(S) \: = \: 
\Gamma (S)^{+} \backslash \Omega _{\widetilde{NS}(S)}^{+}, \: \: \: \: \: 
{\rm where\/} \: \: \: \: 
\Gamma (S) := \Gamma (S, 0). 
\end{equation}

Let 
$l \in I^{d}(\widetilde{NS}(S, B))$ 
be a primitive isotropic vector. 
By using the surjectivity of the period map,   
we shall construct a twisted $K3$ surface 
$(S_{l}, \alpha _{l})$ 
from the cusp 
$[l] \in \mathcal{C}^{d}(S, \alpha )$ 
as follows. 
Let 
$\lambda : (D_{\widetilde{NS}(S, B)}, q) \simeq (D_{T(S, B)}, -q)$ 
be the natural isometry. 
Firstly 
we consider the extended even lattices 
\begin{equation*}\label{extended lattices}
\widetilde{M}_{l} := \Bigl\langle \widetilde{NS}(S, B), \: \frac{l}{d} \Bigr\rangle , \: \: \: \: 
T_{l} := \Bigl\langle T(S, B), \: \lambda (\frac{l}{d}) \Bigr\rangle ,  
\end{equation*}
and a homomorphism 
\begin{equation*}
\alpha _{l} : 
T_{l} 
\twoheadrightarrow 
T_{l}/T(S, B) 
\simeq 
\Bigl\langle \lambda (\frac{l}{d}) \Bigr\rangle 
\simeq  
{\Z}/d{\Z}, \: \: \: \: \: 
\lambda (\frac{l}{d}) \mapsto \bar{1}. 
\end{equation*}
The lattice $\widetilde{M}_{l}$ 
has the orientation induced from that of 
$\widetilde{NS}(S, B)$. 
Since 
$\frac{l}{d} \in I^{1}(\widetilde{M}_{l})$, 
there is an embedding 
$\varphi : U \hookrightarrow \widetilde{M}_{l}$ 
satisfying 
$\varphi (f)=\frac{l}{d}$.  
The orthogonal complement 
$M_{\varphi} := \varphi (U)^{\perp} \cap \widetilde{M}_{l} \; $ 
is of 
${\rm sign\/}(M_{\varphi}) = (1, \rho(S)-1)$. 
We choose the connected component 
$M_{\varphi}^{+}$ 
of the open set 
$\{ v \in (M_{\varphi})_{{\R}} \; | \; (v, v)>0 \}$ 
so that 
${\R}\varphi (e+f) \oplus {\R}v$ 
is of positive orientation 
for all $v \in M_{\varphi}^{+}$. 
From the isometry 
$\bar{\lambda } : (D_{\widetilde{M}_{l}}, q) \simeq (D_{T_{l}}, -q)$, 
we obtain an embedding 
\begin{equation*}
\widetilde{M}_{l} \oplus T_{l} \hookrightarrow \widetilde{\Lambda }_{K3} 
\end{equation*} 
with both 
$\widetilde{M}_{l}$ and $T_{l}$ 
embedded primitively. 
Then 
the orthogonal complement 
$\Lambda _{\varphi} := \varphi (U)^{\perp} \cap \widetilde{\Lambda }_{K3}$ 
is isometric to the $K3$ lattice $\Lambda _{K3}$ 
and is endowed with a period by the sublattice 
$T_{l} \subset \Lambda _{\varphi}$. 
By the surjectivity of the period map (see \cite{K3}), 
there exist a $K3$ surface $S_{l}$ 
and a Hodge isometry 
$\Phi : H^{2}(S_{l}, {\Z}) \simeq \Lambda _{\varphi}$ 
mapping 
the positive cone of $NS(S_{l})$ 
to 
the cone $M_{\varphi}^{+}$. 
Pulling back the homomorphism $\alpha _{l}$ by $\Phi $, 
we obtain the twisted $K3$ surface $(S_{l}, \alpha _{l})$. 
The following lemma can be proved similarly as Lemma 3.2.1 of \cite{Ma1},  
which treats the case $\alpha =1$.

\begin{lemma}\label{well-defined}
The isomorphism class of 
the twisted $K3$ surface $(S_{l}, \alpha _{l})$ 
is uniquely determined by   
the cusp $[l] \in \mathcal{C}^{d}(S, \alpha ) $. 
\end{lemma}

When $S$ is untwisted,  
the lattice $NS(S_{l})$ is isogenus to an overlattice of $NS(S)$ with cyclic quotient.

Let $(S', \alpha ')$ be a given twisted $K3$ surface 
with a B-field lift $B'$ 
such that 
there is an orientation-preserving Hodge isometry 
$\widetilde{\Phi} : 
\widetilde{H}(S', B', {\Z}) \simeq \widetilde{H}(S, B, {\Z})$. 
Then 
we have 
$(S', \alpha ') \in {\rm FM\/}^{d}(S, \alpha )$ 
by Huybrechts-Stellari's theorem (\cite{H-S2}), 
where 
$d={\rm ord\/}(\alpha ')$.

\begin{proposition}\label{key prop}
In the above situation, 
we define the primitive isotropic vector 
$l \in I^{d}(\widetilde{NS}(S, B))$ by 
$l = \widetilde{\Phi} ((0, 0, -1))$. 
Then we have the isomorphism 
$(S_{l}, \alpha _{l}) \simeq (S', \alpha ')$. 
\end{proposition}

\begin{proof}
It suffices to prove the assertion 
when $B' \in \frac{1}{d}H^{2}(S', {\Z})$. 
Let 
\begin{equation*} 
\bar{\lambda} : 
                D_{\widetilde{M}_{l}} \simeq D_{T_{l}}, \: \: \:  \: 
\lambda ' : 
                D_{\widetilde{NS}(S', B')} \simeq D_{T(S', B')},  \: \: \: \:    
\lambda _{0}' : 
                D_{\widetilde{NS}(S')} \simeq D_{T(S')},  
\end{equation*}
be the natural isomorphisms, 
and set 
\begin{equation*}
\widetilde{M}' := 
                  \Bigl\langle \widetilde{NS}(S', B'), \: (0, 0, -\frac{1}{d}) \Bigr\rangle , \: \: \: 
T' := 
                  \Bigl\langle  T(S', B'), \: \lambda '((0, 0, -\frac{1}{d}))   \Bigr\rangle .
\end{equation*} 
We have a natural isometry 
$\bar{\lambda '} : (D_{\widetilde{M}'}, q) \simeq (D_{T'}, -q)$. 
Then 
$\widetilde{\Phi}$ induces 
an orientation-preserving isometry 
$\widetilde{M}' \simeq \widetilde{M}_{l}$ 
and a Hodge isometry 
$T' \simeq T_{l}$. 
By Lemma \ref{twist and disc form} and the commutative diagram 
$$
\CD 
 T(S')  @>{\rm e\/}^{B'}>> T'/T(S', B') @>\simeq >> \langle \lambda '((0, 0, -\frac{1}{d})) \rangle @>\simeq >> {\Z}/d{\Z} 
 \\
@VV\widetilde{\Phi}\circ {\rm e\/}^{B'}V      @VV\widetilde{\Phi}V    @VV\widetilde{\Phi}V        @|     \\ 
 T_{l}  @>>> T_{l}/T(S, B)      @>\simeq >> \langle \lambda (\frac{l}{d}) \rangle @>\simeq >> {\Z}/d{\Z},  
\endCD
$$      
we obtain 
$(\widetilde{\Phi} \circ {\rm e\/}^{B'})^{\ast}\alpha _{l} = \alpha '$.

On the other hand, 
the following diagram commutes by Proposition \ref{commute}. 
$$
\CD 
D_{\widetilde{NS}(S')} @>{\rm e\/}^{B'} \circ \kappa >> D_{\widetilde{M}'} @>\widetilde{\Phi}>> D_{\widetilde{M}_{l}} \\
@VV\lambda _{0}'V       @ VV\bar{\lambda '}V                        @VV\bar{\lambda }V    \\
D_{T(S')}              @>{\rm e\/}^{B'}>> D_{T'}                          @>\widetilde{\Phi}>> D_{T_{l}} .
\endCD
$$
It follows from Proposition \ref{Nikulin} that  
the isometry 
\begin{equation*}
(\widetilde{\Phi} \circ {\rm e\/}^{B'} \circ \kappa ) \oplus (\widetilde{\Phi} \circ {\rm e\/}^{B'}) : 
\widetilde{NS}(S') \oplus T(S') 
\stackrel{\simeq}{\longrightarrow} 
\widetilde{M}_{l} \oplus T_{l} 
\end{equation*}
extends to 
the orientation-preserving Hodge isometry 
\begin{equation*} 
\widetilde{\Psi} : 
                   \widetilde{H}(S', {\Z}) 
                    \stackrel{\simeq}{\longrightarrow} 
                   \widetilde{\Lambda}_{K3}, \: \: \: \: \: 
(0, 0, -1) \mapsto \frac{l}{d}, \: \:  
(\widetilde{\Psi}|_{T(S')})^{\ast} \alpha _{l} = \alpha '.
\end{equation*}
By considering 
$\widetilde{\Psi} : 
H^{0}(S', {\Z}) + H^{4}(S', {\Z}) 
\hookrightarrow 
\widetilde{\Lambda}_{K3}$, 
we obtain an embedding 
$\varphi : U \hookrightarrow \widetilde{M}_{l} \subset \widetilde{\Lambda }_{K3}$ 
with 
$\varphi (f) = \frac{l}{d}$. 
Then  
we have a Hodge isometry 
$\Psi : H^{2}(S', {\Z}) \simeq \Lambda _{\varphi}$ 
which maps the positive cone of $NS(S')$ 
to $M_{\varphi }^{+}$ 
and satisfies 
$(\Psi |_{T(S')})^{\ast} \alpha _{l} = \alpha '$. 
By Lemma \ref{well-defined}, 
we have 
$(S', \alpha ') \simeq (S_{l}, \alpha _{l})$. 
\end{proof}

Let 
${\rm FM\/}^{d}(S, \alpha )^{+}$ 
be the subset of 
${\rm FM\/}^{d}(S, \alpha )$ 
consisting of those partners 
$(S', \alpha ')$ such that 
there are orientation-preserving Hodge isometries  
$\widetilde{H}(S', B', {\Z}) \simeq \widetilde{H}(S, B, {\Z})$. 
It is conjectured by Huybrechts-Stellari (Conjecture 0.2 of \cite{H-S1}) that 
${\rm FM\/}^{d}(S, \alpha )^{+} = {\rm FM\/}^{d}(S, \alpha )$. 
When ${\rm ord\/}(\alpha ) \leq 2$, 
we actually have 
${\rm FM\/}^{d}(S, \alpha )^{+} = {\rm FM\/}^{d}(S, \alpha )$ 
because  
the twisted Mukai lattice admits an orientation-{\it reversing\/} Hodge isometry 
(see \cite{H-S1}).

\begin{proposition}\label{belong to the subset}
The twisted $K3$ surface $(S_{l}, \alpha _{l})$ 
constructed from the cusp 
$[l] \in \mathcal{C}^{d}(S, \alpha ) $ 
belongs to 
${\rm FM\/}^{d}(S, \alpha )^{+}$. 
\end{proposition}

\begin{proof}
We write $l$ degreewise as $l = (r, m, s)$. 

$Case (i)$. 
Assume that $r<0$. 
We can apply 
Yoshioka's theorem (\cite{Yo}. See also \cite{H-S2}) 
on the existence of moduli space of twisted sheaves with Mukai vector $-l$. 
So 
we can find a twisted $K3$ surface $(M, \beta )$ 
and an orientation-preserving Hodge isometry 
$\widetilde{\Phi} : 
\widetilde{H}(M, B', {\Z}) \simeq \widetilde{H}(S, B, {\Z})$ 
with 
$\widetilde{\Phi}((0, 0, 1)) = -l$, 
where 
$B' \in H^{2}(M, {\Q})$ 
is a B-field lift of $\beta $. 
By Proposition \ref{key prop} we have 
$(S_{l}, \alpha _{l}) \simeq (M, \beta ) \in {\rm FM\/}^{d}(S, \alpha )^{+}$.

$Case (ii)$.  
When $r>0$, 
we have 
$(S_{l}, \alpha _{l}) \simeq (S_{-l}, \alpha _{-l}) \in {\rm FM\/}^{d}(S, \alpha )^{+}$ 
by Lemma \ref{well-defined} and the case $(i)$.

$Case (iii)$.
Finally assume that $r=0$. 
By the following lemma, 
there is a vector 
$l' = (r', m', s') \in \Gamma (S, B)^{+} \cdot l$ 
such that  
$r'<0$. 
Thus we have 
$(S_{l}, \alpha _{l}) \simeq (S_{l'}, \alpha _{l'}) \in {\rm FM\/}^{d}(S, \alpha )^{+}$. 
\end{proof}

\begin{lemma}\label{auxiliary}
For a primitive isotropic vector 
$l_{0} = (0, m_{0}, s_{0}) \in \widetilde{NS}(S, B)$, 
there exists a vector 
$l_{1} = (r_{1}, m_{1}, s_{1}) \in O(\widetilde{NS}(S, B))_{0}^{+} \cdot l_{0}$ 
such that 
$r_{1}<0$. 
\end{lemma}

\begin{proof}
Set 
$M := \{ m \in NS(S) \: | \: (m, B) \in {\Z} \}$. 
Take a positive integer  
$d'$ so that 
$d'(1, B, \frac{1}{2}(B, B))$ 
belongs to 
$\widetilde{NS}(S, B)$. 
We can find a vector 
$v\in NS(S)$ 
satisfying the following conditions : 
\begin{equation*}
(m_{0}, v) < 0, \: \: \: (v, v)=0, \: \: \: v \in d'M. 
\end{equation*}
Define the isometry $\varphi $ 
of $\widetilde{NS}(S, B)_{{\Q}} = {\rm e\/}^{B}(\widetilde{NS}(S)_{{\Q}})$ 
by the followings :  
\begin{eqnarray*}
(1, B, \frac{1}{2}(B, B)) & \mapsto &  (1, B, \frac{1}{2}(B, B)),  \\
(0, 0, 1)                 & \mapsto &  (0, 0, 1) + (0, v, (v, B)),  \\ 
(0, m, (m, B))            & \mapsto &  (0, m, (m, B)) + (m, v)\cdot (1, B, \frac{1}{2}(B, B)), \: \: \: \: 
m \in NS(S).
\end{eqnarray*}
Let 
$M':= {\rm e\/}^{B}(M)$ 
and 
$N' := \langle d'(1, B, \frac{1}{2}(B, B)), (0, 0, 1) \rangle $. 
Then 
$M'\oplus N'$ 
is a finite-index sublattice of 
$\widetilde{NS}(S, B)$. 
One can check that 
$\varphi $ preserves 
$M' \oplus N'$ and acts trivially on $D_{M'\oplus N'}$. 
Hence 
$\varphi $ preserves 
$\widetilde{NS}(S, B)$ and acts trivially on $D_{\widetilde{NS}(S, B)}$. 
The $H^{0}(S, {\Z})$-component of 
$\varphi (l_{0})$ is equal to $(m_{0}, v)<0$. 
\end{proof}

Up to now we confirmed that 
the assignment 
$[l] \mapsto (S_{l}, \alpha _{l})$ 
induces a well-defined map 
\begin{equation*}\label{surjection}
\mathcal{C}^{d}(S, \alpha ) 
\longrightarrow 
{\rm FM\/}^{d}(S, \alpha )^{+},  
\end{equation*} 
which is surjective by Proposition \ref{key prop}.  
We prove the injectivity.

\begin{proposition}
Let 
$l_{1}, l_{2} \in I^{d}(\widetilde{NS}(S, B))$ 
be two primitive isotropic vectors. 
If 
$(S_{l_{1}}, \alpha _{l_{1}}) \simeq (S_{l_{2}}, \alpha _{l_{2}})$, 
then  
$l_{2} \in \Gamma (S, B)^{+}\cdot l_{1}$. 
\end{proposition}

\begin{proof}
From the proof of Proposition \ref{belong to the subset}, 
we see the existence of twisted FM partners 
$(M_{i}, \beta _{i})$, $i=1, 2$, 
and orientation-preserving Hodge isometries 
$\widetilde{H}(M_{i}, B_{i}, {\Z}) \simeq \widetilde{H}(S, B, {\Z})$ 
mapping 
$(0, 0, -1)$ to $l_{i}$. 
By the assumption 
we have 
$(M_{1}, \beta _{1}) \simeq (M_{2}, \beta _{2})$. 
Thus 
there is an effective Hodge isometry 
$\Phi : H^{2}(M_{2}, {\Z}) \simeq H^{2}(M_{1}, {\Z})$ 
such that 
the twistings given by 
$B_{1}$ and $\Phi (B_{2})$ coincide. 
So 
we obtain an orientation-preserving Hodge isometry 
$\widetilde{H}(M_{2}, B_{2}, {\Z}) \simeq \widetilde{H}(M_{1}, B_{1}, {\Z})$ 
mapping 
$(0, 0, 1)$ to $(0, 0, 1)$. 
After all, 
we have an orientation-preserving Hodge isometry 
$\widetilde{H}(S, B, {\Z}) \simeq \widetilde{H}(S, B, {\Z})$ 
mapping 
$l_{1}$ to $l_{2}$. 
\end{proof}

Thus we obtain the following theorem.

\begin{theorem}\label{most general}
For a twisted $K3$ surface $(S, \alpha )$, 
the assignment 
$[l] \mapsto (S_{l}, \alpha _{l})$ 
induces a bijection between 
the sets 
\begin{equation*}
\mathcal{C}^{d}(S, \alpha ) 
\simeq 
{\rm FM\/}^{d}(S, \alpha )^{+}. 
\end{equation*}
In particular, 
for a $K3$ surface $S$  
we have a canonical bijection 
\begin{equation*}
\mathcal{C}^{d}(S) 
\simeq 
{\rm FM\/}^{d}(S). 
\end{equation*}
\end{theorem}

In this way  
the set of $0$-dimensional cusps of 
the K\"ahler moduli  
$\mathcal{K}(S)$  
is identified with the set 
$\mathop{\bigcup}_{d} {\rm FM\/}^{d}(S)$ 
of isomorphism classes of the twisted FM partners of $S$.  
For this correspondence,  
it is essential to distinguish two twisted $K3$ surfaces with a common underlying $K3$ surface 
by their twisting classes.

A formula expressing the number $\# {\rm FM\/}^{d}(S)$ is proved in \cite{Ma2}. 
Summing up those formulae over $d \in {\NN}$,  
we obtain a counting formula 
for the $0$-dimensional cusps of the modular variety  $\mathcal{K}(S)$. 
To exhibit the formula, we prepare some notation.   
Let $I(D_{NS(S)})$ 
be the set of isotropic elements of $(D_{NS(S)}, q)$. 
Each element $x \in I(D_{NS(S)})$ 
gives rise to overlattices 
$M_{x}$, $T_{x}$ of $NS(S)$, $T(S)$ respectively, 
and a homomorphism 
$\alpha _{x} : T_{x} \twoheadrightarrow {\Z}/{\rm ord\/}(\alpha ) {\Z}$ 
with 
${\rm Ker\/}(\alpha _{x}) = T(S)$. 
The lattice $T_{x}$ inherits the period from $T(S)$. 
Let 
\begin{equation*}
O_{Hodge}(T_{x}, \alpha _{x}) := 
         \{ g \in O_{Hodge}(T_{x}), \: g^{\ast}\alpha _{x} = \alpha _{x} \} . 
\end{equation*}
For an even lattice $L$, 
the genus of $L$ is denoted by $\mathcal{G}(L)$. 
We define 
\begin{equation*}
\mathcal{G}_{1}(L) := \{ L'\in \mathcal{G}(L) \: |\: O(L')_{0}^{+} \ne O(L')_{0} \} , 
\: \: \: \: \: 
\mathcal{G}_{2}(L) := \mathcal{G}(L) - \mathcal{G}_{1}(L).
\end{equation*}
For a natural number $d \in {\NN}$, 
let 
$\varepsilon (d)=1$ if $d\leq 2$, 
and 
$\varepsilon (d)=2$ if $d\geq 3$.

\begin{theorem}\label{number of cusps}
For a pair 
$(x, M)$ such that 
$x \in I(D_{NS(S)})$ and 
$M \in \mathcal{G}(M_{x})$, 
we write  
\begin{equation*}
\tau (x, M) := \# 
\left( 
O_{Hodge}(T_{x}, \alpha _{x}) \backslash O(D_{M_{x}}) /O(M) 
\right) . 
\end{equation*}
Then 
\begin{equation*}
\#  \: \mathop{\bigcup}_{d} \mathcal{C}^{d}(S)  
\: = \: 
\mathop{\sum}_{x} \Bigl\{ 
\mathop{\sum}_{M} \tau (x, M) + \varepsilon ({\rm ord\/}(x)) \mathop{\sum}_{M'} \tau (x, M')  \Bigr\} . 
\end{equation*}
Here 
$x \in O_{Hodge}(T(S)) \backslash I(D_{NS(S)})$, 
$M \in \mathcal{G}_{1}(M_{x})$, 
and 
$M' \in \mathcal{G}_{2}(M_{x})$. 
\end{theorem}

This counting formula for the cusps is an extension of a formula of Scattone \cite{Sc},     
which in our situation  
is the formula for those $K3$ surfaces whose N\'eron-Severi lattices contain the hyperbolic plane $U$.  
It is the appearance of the subtle arithmetic invariants  $\mathcal{G}_{i}(M_{x})$ and  $\tau (x, M)$ 
that our generalization brings.

\section{Twisted Fourier-Mukai partners and relative Jacobians}\label{sec4}

\subsection{Relative Jacobian}\label{sec4.1}

Let $f : S \to {\proj}^{1}$ 
be an elliptic fibration on a $K3$ surface $S$ 
which does not necessarily admit a section, 
and let $l \in NS(S)$ be the class of its fibres. 
The vector  
$v := (0, l, 0)$  
is a primitive isotropic vector in $\widetilde{NS}(S)$. 
Let 
$J(S/{\proj}^{1})$ 
be the coarse moduli space of stable (with respect to a generic polarization) sheaves on $S$ 
with Mukai vector $v$. 
By a theorem of Yoshioka \cite{Yo}, 
$J(S/{\proj}^{1})$ 
is a $K3$ surface. 
Let $\beta \in {\rm Br\/}(J(S/{\proj}^{1}))$ 
be the obstruction to the existence of a universal sheaf (\cite{Ca2}).   
The order of $\beta $ is equal to ${\rm div\/}(l)$.

\begin{definition}\label{def of relative Jacobian}
We call  $(J(S/{\proj}^{1}), \beta )$ 
the {\it relative Jacobian\/} of the elliptic fibration 
$f : S \to {\proj}^{1}$. 
\end{definition}

The underlying $K3$ surface 
$J(S/{\proj}^{1})$ naturally has the structure of an elliptic fibration 
$\widehat{f} : J(S/{\proj}^{1}) \to {\proj}^{1}$,  
which admits a natural section.  
If $s \in {\proj}^{1}$ is the base locus of a smooth fibre of $f$, 
then 
$\widehat{f}^{-1}(s)$ is identified with ${\rm Pic\/}^{0}(f^{-1}(s))$.  
When the original fibration $f$ admits a section,  
which is exactly the case $\beta =1$,  
we have 
$J(S/{\proj}^{1})/{\proj}^{1} \simeq S/{\proj}^{1}$.

By a theorem of C\u ald\u araru \cite{Ca2}, 
a $1\boxtimes \beta ^{-1}$-twisted universal sheaf on 
$S \times J(S/{\proj}^{1})$  induces an equivalence 
$D^{b}(J(S/{\proj}^{1}), \beta ) \simeq D^{b}(S)$.  
If we denote by $B$ a B-field lift of $\beta $, 
this derived equivalence induces an orientation-preserving Hodge isometry 
$\widetilde{H}(J(S/{\proj}^{1}), B, {\Z}) \simeq \widetilde{H}(S, {\Z})$ 
mapping $(0, 0, 1)$ to $v=(0, l, 0)$.

\begin{theorem}\label{twisted partners and relative Jacobians} 
Let $S$ be a $K3$ surface 
such that $NS(S)$ admits an embedding of the hyperbolic plane $U$.  
Then for each twisted Fourier-Mukai partner $(S', \alpha ')$ of $S$ 
there exists an elliptic fibration 
$f : S \to {\proj}^{1}$ such that 
$(S', \alpha ')$ 
is isomorphic to 
the relative Jacobian $(J(S/{\proj}^{1}), \beta )$ of $f$. 
\end{theorem}

\begin{proof}
By Theorem \ref{most general}, 
there exists a primitive isotropic vector 
$l \in I^{d}(\widetilde{NS}(S))$ 
such that 
$(S', \alpha ') \simeq (S_{l}, \alpha _{l})$, 
where $d = {\rm ord\/}(\alpha ')$. 
Let 
\begin{equation*}
I^{d}(D_{\widetilde{NS}(S)}) := \{  x \in D_{\widetilde{NS}(S)} \: |  \: q(x) \equiv 0, \: {\rm ord\/}(x)=d \: \}. 
\end{equation*} 
Since 
$\widetilde{NS}(S)$ contains $U\oplus U$, 
it follows from Proposition 4.1.3 of \cite{Sc} that 
the map 
\begin{equation}\label{Scattone  bijective}
\Gamma (S)^{+} \backslash I^{d}(\widetilde{NS}(S)) 
\longrightarrow 
r(\Gamma (S)^{+}) \backslash I^{d}(D_{\widetilde{NS}(S)}) ,  \: \: \: \: 
l \mapsto \frac{l}{d}
\end{equation}
is bijective. 
On the other hand, 
the map 
\begin{equation*}\label{Scattone  surjective}
I^{d}(NS(S)) 
\longrightarrow 
I^{d}(D_{NS(S)}) 
\simeq 
I^{d}(D_{\widetilde{NS}(S)}),  \: \: \: \: 
l \mapsto \frac{l}{d}
\end{equation*}
is surjective 
by Lemma 4.1.1 of \cite{Sc}. 
Thus 
there is a primitive isotropic vector 
$l' \in NS(S) \cap (\Gamma (S)^{+}\cdot l)$. 
We have 
$(S', \alpha ') \simeq (S_{l'}, \alpha _{l'})$.

Next,   
we transform $l'$ 
to a nef primitive isotropic vector $l'' \in NS(S)$
by the actions of 
$\{ \pm {\rm id\/} \}$ and 
the reflections with respect to 
$(-2)$-curves on $S$, as in \cite{PS-S}. 
Since 
$l''$ is the fibre class of an elliptic fibration $f : S \to {\proj}^{1}$,  
the twisted FM partner  
$(S_{l''}, \alpha _{l''})$ is isomorphic to the relative Jacobian $(J(S/{\proj}^{1}), \beta )$ of $f$.  
Hence  
we have 
$(S', \alpha ') \simeq (J(S/{\proj}^{1}), \beta )$. 
\end{proof}

\begin{remark}
The bijectivity of the map $(\ref{Scattone  bijective})$ 
can also be deduced from its surjectivity,  which can be proved in an elementary way,    
and the formula for $\# \mathcal{C}^{d}(S)$ (see \cite{Ma2} Corollary 4.4). 
\end{remark}

The above proof might be regarded as an extension of the argument for elliptic fibrations in \cite{PS-S}.  
Theorem \ref{twisted partners and relative Jacobians} is a classification theorem 
for twisted FM partners of $S$.   
The lattice $NS(S)$ contains $U$ 
if and only if 
$S$ admits an elliptic fibration with a section.   
When 
$\rho (S)  \geq 13$,  
we can always embed $U$ into $NS(S)$ by Corollary 1.13.5 of \cite{Ni}.   
Certain $2$-elementary $K3$ surfaces 
give other series of examples satisfying the assumption.

What happens if the assumption that $NS(S)$ contains $U$ is weakened?  
To be precise, 
let $S$ be a $K3$ surface admitting at least one elliptic fibration,    
and $(S', \alpha ')$ be an arbitrary twisted FM partner of $S$.  
Is $(S', \alpha ')$ 
isomorphic to the relative Jacobian of some elliptic fibration on $S$?  
More generally,  
is $(S', \alpha ')$ isomorphic to the twisted moduli space of relative sheaves associated to 
a primitive Mukai vector $(0, \alpha l, \beta )$?   
Here $l$ is the fibre class of some elliptic fibration on $S$,  
$\alpha $ is a natural number,  
and $\beta $ is an integer coprime to $\alpha $.   
The following is a negative example to this question.

\begin{example}\label{partner is not always relative moduli}
Let $S$ be a $K3$ surface such that $NS(S) \simeq U(d) = 
         (
           \begin{smallmatrix} 
              0 & d \\
              d & 0 
            \end{smallmatrix}
          )$  
and $O_{Hodge}(T(S)) = \{ \pm {\rm id \/} \} $.             
The surface $S$ admits exactly two elliptic fibrations,  
whose fibre classes are denoted by $l$ and $m$ respectively.  
For a pair $(\alpha , \beta )$ of coprime integers, 
consider the primitive isotropic vectors in $\widetilde{NS}(S)$ defined by 
\begin{equation*}
v_{\alpha , \beta} = (0, \alpha l, \beta ),  \; \; \; 
v_{\alpha , \beta}' = (0, \alpha m, \beta ).   
\end{equation*}
We put the assumption 
that the twisted FM partner associated to $v_{\alpha , \beta}$ is untwisted,  
which is exactly the case that $\beta $ is coprime to $d$.   
Take integers $\gamma $, $\delta $ satisfying 
$\alpha \gamma d + \beta \delta =1$ 
and define the isotropic vectors in $\widetilde{NS}(S)$  by 
\begin{equation*}
u_{\alpha , \beta}  =    (-\delta ,      \gamma m,          0),  \; \; \; 
l_{\alpha , \beta}  =    (0,               \delta l,   -\gamma d),    \; \; \; 
m_{\alpha , \beta}   =    (\alpha d,  \beta m,               0).  
\end{equation*}   
One checks that 
\begin{equation*}
\langle  v_{\alpha , \beta}, u_{\alpha , \beta} \rangle \simeq U,  \; \; \; 
\langle  l_{\alpha , \beta}, m_{\alpha , \beta} \rangle \simeq U(d),  \; \; \; 
\langle  v_{\alpha , \beta}, u_{\alpha , \beta} \rangle \perp \langle  l_{\alpha , \beta}, m_{\alpha , \beta} \rangle .    
\end{equation*} 
Therefore we obtain an isometry $\varphi _{\alpha , \beta} \in O(\widetilde{NS}(S))$
which satisfies 
$\varphi _{\alpha , \beta}(v_{\alpha , \beta}) =  (0, 0, 1)$  
and acts on the discriminant group 
$D_{NS(S)} = \langle \frac{l}{d}, \frac{m}{d} \rangle $ by  
$(\frac{l}{d},  \frac{m}{d}) \mapsto (\beta \frac{l}{d},  \beta ^{-1}\frac{m}{d})$.     
Since any isometry of $\widetilde{NS}(S)$ fixing $(0, 0, 1)$ 
must preserve the subset 
$\langle \frac{l}{d} \rangle \cup \langle \frac{m}{d} \rangle $  
in $D_{NS(S)}$,   
we see that any isometry 
$\varphi \in O(\widetilde{NS}(S))$ satisfying 
$\varphi (v_{\alpha , \beta}) = (0, 0, 1)$  
must preserve the subset 
$\langle \frac{l}{d} \rangle \cup \langle \frac{m}{d} \rangle $  
in $D_{NS(S)}$.  
By symmetry, 
the same holds for the vector $v_{\alpha , \beta}'$.

When  $d$ can be divided by at least two primes,  
there exists an isometry $\varphi \in O(\widetilde{NS}(S))$ 
which does not preserve the subset 
$\langle \frac{l}{d} \rangle \cup \langle \frac{m}{d} \rangle $  
in $D_{NS(S)}$.   
If we put $v := \varphi ((0, 0, 1))$,  
then $v$ is not $\Gamma (S)^{+}$-equivalent to 
$ v_{\alpha , \beta}$ nor  $v_{\alpha , \beta}'$  
for any $(\alpha , \beta )$.  
In other words, 
the FM partner associated to $v$ cannot be 
realized as moduli space of relative sheaves.   
\end{example}

Bridgeland and Maciocia (\cite{B-M} Proposition 4.4) proved that     
for a minimal surface $S$ of {\it non-zero\/} Kodaira dimension which admits an elliptic fibration $f : S \to C$,  
every untwisted FM partner of $S$ is isomorphic to a moduli space of relative sheaves supported on the fibers of $f$.     
Example \ref{partner is not always relative moduli}  
shows that the analogous statement does no longer holds 
for a general elliptic $K3$ surface, 
even if one allows all elliptic fibrations on the surface.

\begin{remark}
If we do not restrict to untwisted $K3$ surfaces,  
we can give negative examples to the above question more handily  
by analyzing the surjective map 
$I^{d}(\widetilde{NS}(S)) \to I^{d}(D_{NS(S)})$, $d>1$.
\end{remark}

\subsection{Applications to rational maps}\label{sec4.2}

Among the twisted FM partners of a $K3$ surface $S$,   
the relative Jacobians are related rather directly to the geometry of $S$.   
We give an application of Theorem \ref{twisted partners and relative Jacobians} 
to rational maps between $K3$ surfaces.

\begin{lemma}\label{rational map to twisted FM partner}
Let $S$ be a $K3$ surface and 
assume that $NS(S)$ contains $U$. 
Then every twisted Fourier-Mukai partner 
$(S' , \alpha ')$ of $S$ 
admits a rational map $S \dashrightarrow  S'$ 
of degree ${\rm ord\/}(\alpha ')^{2}$. 
\end{lemma}

\begin{proof}
By Theorem \ref{twisted partners and relative Jacobians} 
there exists an elliptic fibration $f : S \to {\proj}^{1}$ such that 
$S'$ is isomorphic to the $K3$ surface $J(S/{\proj}^{1})$ underlying the relative Jacobian of $f$.  
Let $l \in NS(S)$ be the fibre class of $f$.  
The order of $\alpha '$ is equal to the divisibility $d$ of $l$ in $NS(S)$.  
It suffices to construct a rational map 
$S \dashrightarrow J(S/{\proj}^{1})$ 
of degree $d^{2}$.  
Take a line bundle $M$ on $S$ such that $(M.l)=d$. 
Let $U \subset {\proj}^{1}$ be a Zariski open set 
such that $f : f^{-1}(U) \to U$ is smooth.  
We can define a morphism 
$f^{-1}(U) \to J(S/{\proj}^{1})$ 
by setting 
\begin{equation}\label{relative multiplication map}
P \mapsto \mathcal{O}_{f^{-1}(f(P))}(dP)\otimes M^{-1}, \: \: \: \: P \in f^{-1}(U), 
\end{equation}
the degree of which is obviously equal to $d^{2}$. 
\end{proof}

We provide an example of the rational map given in Lemma \ref{rational map to twisted FM partner}.

\begin{example}\label{example of multiplication map} 
Let $g : A \to E$ be an elliptic fibration on an Abelian surface $A$ where $E$ is an elliptic curve,  
and $f : {\rm Km\/} (A) \to {\proj}^{1}$ be the associated elliptic fibration on the Kummer surface ${\rm Km\/} (A)$.   
We take a line bundle $M$ on ${\rm Km\/} (A)$ such that 
the degree of $M$ on an $f$-fiber is equal to the divisibility $d$ of the fiber class of $f$ in $NS({\rm Km\/} (A))$.  
The pullback of $M$ by the rational quotient map $A \dashrightarrow {\rm Km\/} (A)$ extends to a line bundle $L$ on $A$.   
The degree of $L$ on a $g$-fiber is also equal to $d$.  
Let $\widehat{g} : B \to E$ be the elliptic surface associated to the relative Jacobian of $g$.  
The surface $B$ is isomorphic to the product of $E$ and a $g$-fiber.  
We can use the line bundle $L$ to define a morphism 
$\varphi : A \to B$ by $P \mapsto \mathcal{O}_{g^{-1}(g(P))}(dP)\otimes L^{-1}$,  
where $P \in A$.  
We choose the identity point of $B$ so that $\varphi $ is a homomorphism between Abelian surfaces. 
Then the elliptic Kummer surface  $\widehat{f} : {\rm Km\/} (B) \to {\proj}^{1}$ induced from $\widehat{g}$  
is canonically isomorphic to the elliptic $K3$ surface  associated to the relative Jacobian of $f$.  
Under this identification,  
the rational map $\psi : {\rm Km\/} (A) \dashrightarrow {\rm Km\/} (B)$ induced from the homomorphism $\varphi $ 
is given by the correspondence $(\ref{relative multiplication map})$ in Lemma \ref{rational map to twisted FM partner}.

We describe the elimination of the indeterminacy of $\psi $.  
Let $A_{2}$, $B_{2}$ be the sets of $2$-division points of $A$, $B$ respectively, 
and $\widetilde{A}$ be the blow-up of $A$ at $\varphi ^{-1}(B_{2})$.  
The set $T$ consists of $16d^{2}$ $2d$-division points and contains $A_{2}$.  
The inverse morphism $-1_{A}$ of $A$ extends to 
an involution $\iota $ on $\widetilde{A}$.  
We denote by $X$ the quotient surface $\widetilde{A}/\langle \iota \rangle $.   
The natural morphism $\pi : X \to {\rm Km\/} (A)$ is the blow-up at the point set $(\varphi ^{-1}(B_{2})\backslash A_{2})/\langle -1_{A} \rangle $,   
and the finite morphism $\widetilde{\psi} : X \to {\rm Km\/} (B)$  induced from $\varphi $ 
is the elimination of the indeterminacy of $\psi $.   
The ramification divisor of $\widetilde{\psi}$ is the exceptional divisor of $\pi $,   
namely the $8(d^{2}-1)$ $(-1)$-curves.   
The ramification index of $\widetilde{\psi}$ at each $(-1)$-curve is $2$.   
Let 
$\{ 2E_{i}^{A} +E_{i1}^{A}+E_{i2}^{A}+E_{i3}^{A}+E_{i4}^{A} \} _{i=1,\cdots ,4} $  
and 
$\{ 2E_{i}^{B} +E_{i1}^{B}+E_{i2}^{B}+E_{i3}^{B}+E_{i4}^{B} \} _{i=1,\cdots ,4} $ 
be the four $I_{0}^{\ast}$ singular fibers of $f$, $\widehat{f}$ respectively.    
All components $E_{i}^{\ast}$,  $E_{ij}^{\ast}$ are $(-2)$-curves, 
and the components $E_{i}^{A}$, $E_{i}^{B}$ are the quotients of fibers of $g$, $\widehat{g}$ respectively.   
The surface $X$ is obtained by blowing up $2(d^{2}-1)$ points on each curve $E_{i}^{A}$, $1\leq i\leq 4$.   
The branch curve of $\widetilde{\psi}$ is contained in the curve $\bigcup _{i,j} E_{ij}^{B}$.   
To observe the branching further,   
we consider the case of odd $d$ and the case of even $d$ separately.   
Renumbering the $i$ if necessary, 
we may assume that the image of $E_{i}^{A}$ by $\psi $ is $E_{i}^{B}$.  
When $d$ is odd, 
the curve $\bigcup _{j=1}^{4} E_{ij}^{A}$ is mapped by $\psi $ isomorphically to $\bigcup _{j=1}^{4} E_{ij}^{B}$.   
We may assume that $\psi (E_{ij}^{A}) = E_{ij}^{B}$.  
Then for each $j$ 
the inverse image $\widetilde{\psi}^{-1}(E_{ij}^{B})$ consists of 
the proper transform of $E_{ij}^{A}$ and $\frac{1}{2}(d^{2}-1)$ ramifying $(-1)$-curves.   
When $d$ is even,   
after renumbering the $j$ for $E_{ij}^{B}$,   
we have 
$\psi (E_{ij}^{A}) = E_{i1}^{B}$ for all $j$.  
Then
$\widetilde{\psi}^{-1}(E_{i1}^{B})$ consists of 
the proper transform of  $\bigcup _{j=1}^{4} E_{ij}^{A}$  and  $\frac{1}{2}(d^{2}-4)$ ramifying $(-1)$-curves. 
On the other hand,  
for $2\leq j\leq 4$, $\widetilde{\psi}^{-1}(E_{ij}^{B})$ consists of $\frac{1}{2}d^{2}$ ramifying $(-1)$-curves.  
We remark that when $d=2$, 
$\psi $ is the Galois covering
for the symplectic action of the group ${\rm Ker\/}(\varphi ) \simeq ({\Z}/2{\Z})^{\oplus 2}$ on ${\rm Km\/} (A)$.  
\end{example}

Following Mukai's approach in \cite{Mu},   
we deduce the next proposition.

\begin{proposition}\label{Hodge embedding and rational map}
Let $S$ and $S'$ be $K3$ surfaces with 
$\rho (S) = \rho (S')$. 
Assume that $NS(S)$ contains $U$. 
If there is a Hodge embedding $T(S) \hookrightarrow T(S')$, 
then there exists a rational map $S \dashrightarrow S'$ 
of degree $[ T(S') : T(S) ]^{2}$. 
\end{proposition}

\begin{proof}
We regard $T(S)$ as a finite-index sublattice of $T(S')$ 
via the Hodge embedding. 
By the invariant factor theorem, 
there exists a filtration 
\begin{equation*} 
T(S) = T_{0} \subset T_{1} \subset \cdots \subset T_{n} = T(S')  
\end{equation*}
such that 
$T_{i}/T_{i-1}$ are cyclic. 
By using similar arguments as in 
the construction of $S_{l}$ in Lemma \ref{well-defined}, 
we can find for each $i$ a $K3$ surface $S_{i}$ 
such that $T(S_{i})$ is Hodge isometric to $T_{i}$  
and that $NS(S_{i})$ is isometric to an overlattice of $NS(S)$.    
By the assumption on $NS(S)$,    
the lattice $NS(S_{i})$ admits an embedding of $U$. 
Since there exist 
Hodge embeddings $T(S_{i-1}) \hookrightarrow T(S_{i})$ 
of cyclic quotients, 
we can find twistings 
$\alpha _{i} \in {\rm Br\/}(S_{i})$ 
such that 
$T(S_{i}, \alpha _{i})$ 
are Hodge isometric to 
$T(S_{i-1})$.
By \cite{H-S2}, 
we have 
$D^{b}(S_{i}, \alpha _{i}) \simeq D^{b}(S_{i-1})$. 
Now the assertion follows from 
Lemma \ref{rational map to twisted FM partner}. 
\end{proof}

Two $K3$ surfaces $S_{+}$ and $S_{-}$ are 
{\it isogenous\/} 
(in the sense of Mukai \cite{Mu}) 
if there exists an algebraic cycle 
$Z \in H^{4}(S_{+} \times S_{-}, {\Q})$ 
such that the correspondence 
\begin{equation*}\label{correspondence} 
\Phi _{Z} : H^{2}(S_{+}, {\Q}) \to H^{2}(S_{-}, {\Q}), \: \: \: \:  
l \mapsto (\pi _{-})_{\ast} (Z \wedge \pi _{+}^{\ast}l), 
\end{equation*}
is a Hodge isometry.

\begin{proposition}\label{isogenous}
Let $S_{+}$ and $S_{-}$ be $K3$ surfaces with 
$\rho (S_{\pm}) \geq 13$. 
The following two conditions are equivalent: 

$(1)$ 
$S_{+}$ and $S_{-}$ are isogenous. 

$(2)$ 
There exist a $K3$ surface $S_{0}$ and 
rational maps $S_{0} \dashrightarrow S_{+}$,  $S_{0} \dashrightarrow S_{-}$ of square degrees.
\end{proposition}

\begin{proof} 
By Mukai's theorem (Corollary 1.10 of \cite{Mu}), 
$S_{+}$ and $S_{-}$ are isogenous 
if and only if 
$T(S_{+})_{{\Q}}$ and $T(S_{-})_{{\Q}}$ are Hodge isometric.

$(1) \Rightarrow (2)$ : 
We identify $T(S_{+})_{{\Q}}$ and $T(S_{-})_{{\Q}}$ by a Hodge isometry. 
Let $T_{0} := T(S_{+}) \cap T(S_{-})$, 
which is of finite index in both $T(S_{+})$ and $T(S_{-})$. 
The lattice $T_{0}$ is endowed with the period. 
Since ${\rm rk\/}(T_{0}) \leq 9$, 
the lattice $T_{0}$ can be embedded primitively into 
the $K3$ lattice $\Lambda _{K3}$ by \cite{Ni}. 
By the surjectivity of the period map, 
there exist a $K3$ surface $S_{0}$ and 
a Hodge isometry $T(S_{0}) \simeq T_{0}$. 
The lattice $NS(S_{0})$ admits an embedding of $U$, 
because $\rho (S_{0}) \geq 13$. 
Hence the claim follows from Proposition \ref{Hodge embedding and rational map}.

$(2) \Rightarrow (1)$ : 
Let $f_{\pm} : S_{0} \dashrightarrow S_{\pm}$ 
be rational maps of degree $d_{\pm}^{2}$. 
We have Hodge embeddings 
$f_{\pm}^{\ast} : T(S_{\pm})(d_{\pm}^{2}) \hookrightarrow T(S_{0})$ 
of finite indices. 
Since 
the lattices $T(S_{\pm})(d_{\pm}^{2})$ 
can be embedded into 
the lattices $T(S_{\pm})$ 
by the multiplications by $d_{\pm}$, 
it follows that 
$T(S_{\pm})_{{\Q}}$ are Hodge isometric to $T(S_{0})_{{\Q}}$ .
\end{proof}

By the proof, 
the $K3$ surfaces $S_{+}$ and $S_{-}$ 
are obtained from $S_{0}$ 
by taking certain relative Jacobians successively.

\begin{remark}
Inose \cite{In} introduced another notion of isogeny for singular $K3$ surfaces, 
i.e. $K3$ surfaces with Picard number $20$.  
Two singular $K3$ surfaces $S_{+}$ and $S_{-}$ are defined to be isogenous in the sense of Inose 
if one of the following three equivalent conditions is satisfied:    
\begin{enumerate}
 \item[(1)]    There exists a dominant rational map  $S_{+} \dashrightarrow S_{-}$.
 \item[(2)]    There exists a dominant rational map  $S_{-} \dashrightarrow S_{+}$.      
 \item[(3)]   There exists a Hodge isometry $T(S_{+})(n)_{{\mathbb Q}} \simeq  T(S_{-})_{{\mathbb Q}}$ for some natural number $n$.  
\end{enumerate} 
For singular $K3$ surfaces (and also for $K3$ surfaces with Picard number $19$) 
Inose's notion of isogeny contains that of Mukai 
and is a direct analogue of the notion of isogeny for Abelian varieties.  
Unfortunately,  the author does not know successful extension of Inose's notion of isogeny 
to $K3$ surfaces with Picard number  $\leq 16$.  
See also \cite{Ni2}. 
\end{remark}

\section{The case of Picard number $1$}\label{sec5}

Let $S$ be a $K3$ surface with 
$NS(S)={\Z}H$, 
$(H, H)=2n>0$. 
In this section 
we shall calculate for $S$ 
the correspondence between 
the twisted FM partners 
and the $0$-dimensional cusps concretely.    
The set 
${\rm FM\/}^{d}(S)$ 
is calculated in \cite{Ma2} 
as a set of moduli spaces of sheaves on $S$ 
with explicit Mukai vectors, 
twisted by the obstruction classes. 
On the other hand, 
Dolgachev showed in \cite{Do} Theorem 7.1 that 
the group $O(\widetilde{NS}(S))_{0}^{+}$ is isomorphic to the Fricke modular group, 
by direct calculations based on an isomorphism 
$PSO(1, 2) \simeq PSL_{2}({\R})$. 
Hence the K\"ahler moduli $\mathcal{K} (S)$ is isomorphic to the Fricke modular curve.  
We shall describe the groups $\Gamma (S)$ and $O(\widetilde{NS}(S))$   
via quaternion orders, which are nicely compatible with the discriminant form.  
The tube domain realization of $\Omega _{\widetilde{NS}(S)}^{+}$ 
induces an explicit isomorphism between the modular curves.

\subsection{Fricke modular curves}\label{sec5.1}

Let 
\begin{equation*}\label{the order}
\mathcal{O} :=  
\left\{  \left. 
\begin{pmatrix}
  a       & 2b  \\ 
 2nc   & d   
\end{pmatrix} 
\in M_{2}({\Q}) 
\: \: \right|  \: \: 
a, b, c, d \in {\Z}, \: \: a+d \in 2{\Z}
\right\} 
\end{equation*}
and  
$
\mathcal{O}_{0} := \{  A\in \mathcal{O}, \: {\rm Tr\/}(A)=0 \}  
$. 
A natural decomposition 
$
\mathcal{O} = 
{\Z}
(
\begin{smallmatrix}
1 & 0 \\
0 & 1 
\end{smallmatrix}
)
\oplus 
\mathcal{O}_{0} 
$ 
holds, 
and $\mathcal{O}$ has a basis 
\begin{equation}\label{basis}
\left\{ 
\begin{pmatrix}
1 & 0 \\
0 & 1 
\end{pmatrix} , \: 
\begin{pmatrix}
-1 & 0 \\
0 & 1 
\end{pmatrix} , \: 
\begin{pmatrix}
0 & 0 \\
-2n & 0 
\end{pmatrix} , \: 
\begin{pmatrix}
0 & 2 \\
0 & 0 
\end{pmatrix} \: 
\right\} .
\end{equation} 
The ${\Z}$-module $\mathcal{O}$ 
is of rank $4$ 
and 
is closed under multiplication. 
In other words, 
$\mathcal{O}$ is an {\it order\/}.

We define a quadratic form on 
the ${\Q}$-vector space 
$V :=  \{  A\in M_{2}({\Q}), \: {\rm Tr\/}(A)=0 \}$ 
by 
$(A, B) := -\frac{1}{2}{\rm Tr\/}(AB)$,  where $A, B \in V$.

\begin{lemma}
The lattice $\mathcal{O}_{0} \subset V$ 
is isometric to 
$\Lambda (-2n)$, 
where 
$\Lambda = \widetilde{NS}(S)^{\vee}$. 
\end{lemma}

\begin{proof}
In fact, 
an isometry is given by the assignment 
\begin{equation}\label{isometry between order and NS}
a
\begin{pmatrix}
0 & 0 \\
-2n & 0 
\end{pmatrix} 
+ b  
\begin{pmatrix}
-1 & 0 \\
0 & 1 
\end{pmatrix} 
+ c 
\begin{pmatrix}
0 & 2 \\
0 & 0 
\end{pmatrix} 
\: \: 
\mapsto 
(a, \; b\frac{H}{2n}, \; c) \: \in \widetilde{NS}(S)^{\vee}.
\end{equation}
\end{proof}

For a matrix 
$\gamma \in GL_{2}({\Q})$ 
the adjoint action ${\rm Ad\/}_{\gamma}$ on $V$ defined by 
$A \mapsto \gamma A \gamma ^{-1}, \:  A \in V$,
is an isometry of $V$. 
Then 
we obtain the isomorphism 
\begin{equation*}
GL_{2}({\Q})/{\Q}^{\times} \simeq SO(V),  \: \: \: \: \: 
\gamma \mapsto {\rm Ad\/}_{\gamma} . 
\end{equation*}
The normalizer $N(\mathcal{O}_{0}) \subset GL_{2}({\Q})$ 
of $\mathcal{O}_{0}$ in $V$ 
coincides with 
the normalizer $N(\mathcal{O}) \subset GL_{2}({\Q})$ 
of $\mathcal{O}$ in $M_{2}({\Q})$.  
Thus 
we obtain the isomorphisms 
\begin{equation}\label{isometry and normalizer}
O(\widetilde{NS}(S)^{\vee}) 
\simeq 
O(\mathcal{O}_{0})  
\simeq 
\langle N(\mathcal{O})/{\Q}^{\times}, \: -{\rm id\/} \rangle . 
\end{equation} 
To calculate the normalizer $N(\mathcal{O})$ 
we consider the following order $\mathcal{O}'$ 
called an  
{\it Eichler order\/}: 
\begin{equation}\label{Eichler order}
\mathcal{O}' :=  
\left\{  \left. 
\begin{pmatrix}
  a       & b  \\ 
 nc   & d   
\end{pmatrix} 
\in M_{2}({\Q}) 
\: \: \right|  \: \: 
a, b, c, d \in {\Z} 
\right\} . 
\end{equation}
The ${\Z}$-module $2\mathcal{O}'$ 
is a submodule of $\mathcal{O}$ of index $2$.



\begin{proposition} 
We have 
$N(\mathcal{O}) = N(\mathcal{O}')$. 
\end{proposition}

\begin{proof}
Let 
$\gamma \in N(\mathcal{O}')$. 
Since 
$
{\rm Ad\/}_{\gamma}
(
\begin{smallmatrix}
1 & 0 \\
0 & 0 
\end{smallmatrix}
), 
{\rm Ad\/}_{\gamma}
(
\begin{smallmatrix}
0 & 1 \\
0 & 0 
\end{smallmatrix}
),
{\rm Ad\/}_{\gamma}
(
\begin{smallmatrix}
0 & 0 \\
n & 0 
\end{smallmatrix}
) 
\in \mathcal{O}'$, 
then we have 
$
{\rm Ad\/}_{\gamma}
(
\begin{smallmatrix}
2 & 0 \\
0 & 0 
\end{smallmatrix}
), 
{\rm Ad\/}_{\gamma}
(
\begin{smallmatrix}
0 & 2 \\
0 & 0 
\end{smallmatrix}
),
{\rm Ad\/}_{\gamma}
(
\begin{smallmatrix}
0 & 0 \\
2n & 0 
\end{smallmatrix}
) 
\in 2\mathcal{O}' \subset \mathcal{O}$.  
Thus   
$N(\mathcal{O}') \subset N(\mathcal{O})$.

Let 
$\gamma \in N(\mathcal{O})$. 
Then 
${\rm Ad\/}_{\gamma}(2\mathcal{O}')$ 
is an index $2$ submodule of $\mathcal{O}$, 
and the ${\Z}$-module 
$\frac{1}{2}{\rm Ad\/}_{\gamma}(2\mathcal{O}') = {\rm Ad\/}_{\gamma}(\mathcal{O}')$ 
is closed under multiplication. 
Among 
$2^{4}-1 = 15$ index $2$ submodules $L$ of $\mathcal{O}$, 
$2\mathcal{O}'$ is characterized by 
the property that 
$\frac{1}{2}L$ is closed under multiplication. 
Hence 
${\rm Ad\/}_{\gamma}(2\mathcal{O}') = 2\mathcal{O}'$ 
and so 
$\gamma \in N(\mathcal{O}')$. 
\end{proof}

We shall describe a structure of the normalizer $N(\mathcal{O}')$ following \cite{Vi} 
and then calculate the action of $N(\mathcal{O}')$ on the discriminant group $D_{\widetilde{NS}(S)}$.   
For $n>1$  
let  
$n=\prod_{i=1}^{\tau (n)}p_{i}^{e_{i}}$ 
be the prime decomposition of $n$, 
where 
$\tau (n)$ is the number of the prime divisors of $n$.  
Each element 
$\sigma \in ({\Z}/2{\Z})^{\tau (n)}$ 
corresponds to 
a pair of coprime natural numbers 
$(N_{\sigma}, M_{\sigma})$ 
satisfying $N_{\sigma}M_{\sigma}=n$.  
Choose integers 
$a_{\sigma}, b_{\sigma}$ 
such that 
$a_{\sigma}N_{\sigma}-b_{\sigma}M_{\sigma} =1$ 
and put 
\begin{equation*}\label{Atkin-Lehner involutions}
\gamma _{\sigma} 
:= 
\begin{pmatrix}
a_{\sigma}N_{\sigma} & b_{\sigma} \\ 
n                                            & N_{\sigma}
\end{pmatrix} .
\end{equation*}
We have 
${\rm det\/}( \gamma _{\sigma} ) = N_{\sigma}$.  
For the $\sigma $ with 
$(N_{\sigma}, M_{\sigma}) = (n, 1)$, 
we take 
$(a_{\sigma}, b_{\sigma}) = (0, -1)$ 
especially. 
Then 
$
\gamma _{\sigma} = 
(
\begin{smallmatrix}
1  & 0 \\
-n & 1 
\end{smallmatrix}
)
(
\begin{smallmatrix}
0  & -1 \\
n  & 0 
\end{smallmatrix}
)
$.

\begin{proposition}[cf.\cite{Vi}]\label{description of the normalizer}
Let $(\mathcal{O}')^{\times}$ be the unit group of $\mathcal{O}'$.
For $n>1$   
we have 
$ N(\mathcal{O}')/{\Q}^{\times}(\mathcal{O}')^{\times} \simeq ({\Z}/2{\Z})^{\tau (n)}$ 
with the quotient group 
$({\Z}/2{\Z})^{\tau (n)}$ 
represented by the set 
$\{ \gamma _{\sigma} \} _{\sigma}$. 
If $n=1$, 
we have 
$
N(\mathcal{O}') = N(M_{2}({\Z})) = 
{\Q}^{\times}\langle SL_{2}({\Z}), 
         (
           \begin{smallmatrix} 
              1 & 0 \\
              0 & -1 
            \end{smallmatrix}
          ) 
          \rangle 
$. 
\end{proposition}

\begin{proof}
Let $n>1$. 
For $i \leq \tau (n)$ we denote 
$g_{i} : = \left( \begin{smallmatrix}
                                                                                  0                 &     1 \\
                                                                              p_{i}^{e_{i}}   &      0  
                                  \end{smallmatrix} \right) $. 
The two maximal orders 
$\mathcal{O}^{+}_{i} := M_{2}({\Z}_{p_{i}})$  
and 
$\mathcal{O}^{-}_{i} := g_{i} \mathcal{O}^{+}_{i}g_{i}^{-1}$ 
in $M_{2}({\Q}_{p_{i}})$ 
satisfy 
$\mathcal{O}'\otimes {\Z}_{p_{i}} = \mathcal{O}^{+}_{i} \cap \mathcal{O}^{-}_{i}$.   
It is known (see \cite{Vi} Chapitre III) that 
we have a surjective homomorphism 
$N(\mathcal{O}') \to ({\Z}/2{\Z})^{\tau (n)}$  
by looking for each $i$ whether the action preserves or exchanges $\mathcal{O}_{i}^{\pm }$,     
and that the kernel is given by ${\Q}^{\times}(\mathcal{O}')^{\times}$.  
Direct calculations show that 
$\gamma _{\sigma}$  
exchanges $\mathcal{O}_{i}^{\pm }$  for those $i$ with $p_{i}|N_{\sigma}$ 
and preserves $\mathcal{O}_{j}^{\pm }$ for those $j$ with $p_{j}|M_{\sigma}$. 
\end{proof}

Let $\Gamma _{0}(n) := \mathcal{O}' \cap SL_{2}({\Z})$ as usual. 
We have 
$(\mathcal{O}')^{\times}  =  
\langle \Gamma _{0}(n), 
(
\begin{smallmatrix} 
1 & 0 \\
0 & -1 
\end{smallmatrix}
) 
\rangle $. 
By $(\ref{isometry and normalizer})$ 
and Proposition \ref{description of the normalizer}
we obtain a description of the full isometry group 
$O(\widetilde{NS}(S)^{\vee}) = O(\widetilde{NS}(S))$.  
To find its subgroup 
\begin{equation*}
\Gamma (S)^{+} = O(\widetilde{NS}(S))_{0}^{+} \times \{ \pm {\rm id\/} \} 
\subset 
 O(\widetilde{NS}(S)^{\vee}), 
 \end{equation*}
 we study the actions of 
 $\Gamma _{0}(n)$, 
 $
 (
 \begin{smallmatrix}
 1 & 0 \\
 0 & -1 
 \end{smallmatrix}
 )
$, 
and  
$\{ \gamma _{\sigma} \}$ 
on the discriminant group $D_{\widetilde{NS}(S)}$. 
By the correspondence $(\ref{isometry between order and NS})$, 
it suffices to observe the actions to 
$(
\begin{smallmatrix}
-1 & 0 \\ 
0  & 1 
\end{smallmatrix}
)$ 
modulo 
$
{\Z}( 
          \begin{smallmatrix} 
                -2n & 0 \\ 
                  0     & 2n 
           \end{smallmatrix}
         ) 
\oplus 
{\Z}( 
           \begin{smallmatrix} 
                   0       & 0 \\ 
                   2n     & 0 
           \end{smallmatrix}
         ) 
\oplus 
{\Z}( 
         \begin{smallmatrix} 
                0     & 2 \\ 
                0     & 0 
           \end{smallmatrix}
         ) 
$. 
Then 
$\Gamma _{0}(n)$ and  
$(
\begin{smallmatrix}
1  & 0 \\ 
0  & -1 
\end{smallmatrix}
)$ 
act trivially on $D_{\widetilde{NS}(S)}$. 
Let $n>1$.   
The isometry $\gamma _{\sigma}$ 
acts on $D_{\widetilde{NS}(S)} \simeq {\Z}/2n{\Z}$ 
as multiplication by 
$a_{\sigma}N_{\sigma} + b_{\sigma}M_{\sigma}$.   
For odd $p_{i}$, 
$\gamma _{\sigma}$ 
acts as 
$ - {\rm id\/} $ (resp. $ {\rm id\/} $)  
on the component 
${\Z}/p_{i}^{e_{i}}{\Z}$ of $D_{\widetilde{NS}(S)}$ 
if $p_{i}|N_{\sigma}$ (resp. $p_{i}|M_{\sigma}$).  
For $p=2$, 
$\gamma _{\sigma}$ 
acts as 
$ - {\rm id\/} $ (resp. $ {\rm id\/} $)  
on the component 
${\Z}/2^{e+1}{\Z}$  
if $2|N_{\sigma}$ (resp. $2|M_{\sigma}$).  
It is worth noting that  
$\gamma _{\sigma}$ acts as $ - {\rm id\/} $ 
on the component ${\Z}/p_{i}^{e_{i}(+1)}{\Z}$ 
if and only if 
$\gamma _{\sigma}$ exchanges the orders $\mathcal{O}_{i}^{\pm }$.

In conclusion,  
we have 
\begin{equation*}
\Gamma (S) 
\simeq 
\Bigl\langle 
              \langle 
                       {\Q}^{\times}\Gamma _{0}(n),  \: 
                        (\begin{smallmatrix}
                                        1 & 0 \\
                                        0 & -1 
                           \end{smallmatrix}), \: 
                         (\begin{smallmatrix}
                                        0 & -1 \\
                                        n &  0 
                           \end{smallmatrix})
               \rangle  / {\Q}^{\times},  \: 
       -{\rm id\/}
\Bigr\rangle  , 
\end{equation*}
and 
therefore 
\begin{eqnarray}\label{final description}
                     \Gamma (S)^{+} 
& \simeq  & 
                     \Bigl\langle 
                                    \langle 
                                           {\Q}^{\times}\Gamma _{0}(n),  \: 
                                                (\begin{smallmatrix}
                                                              0 & -1 \\
                                                               n &  0 
                                                   \end{smallmatrix})
                                     \rangle  / {\Q}^{\times},  \: 
                                    -{\rm id\/}
                      \Bigr\rangle    \nonumber\\
& \simeq  &  
                      \Bigl\langle 
                                      \langle 
                                             \Gamma _{0}(n),  \: 
                                               (\begin{smallmatrix}
                                                               0                  & -\sqrt{n}^{-1} \\
                                                                \sqrt{n}    &  0 
                                                 \end{smallmatrix})
                                      \rangle  ,  \: 
                                 -{\rm id\/} 
                       \Bigr\rangle  .
\end{eqnarray} 
The last isomorphism is induced  by the projection to $SL_{2}({\R})$.  
The group 
$
\Gamma _{0}(n)^{+} := 
\langle  \Gamma _{0}(n),  \:  (\begin{smallmatrix}
                                                               0                  & -\sqrt{n}^{-1} \\
                                                                \sqrt{n}    &  0 
                                                 \end{smallmatrix}) 
\rangle $ 
is called 
the {\it Fricke modular group\/}.

Consider the tube domain realization 
\begin{equation}\label{eqn:tube domain realization}
 {\mathbb H} 
\stackrel{\simeq}{\longrightarrow} 
\Omega _{\widetilde{NS}(S)}^{+}, \: \: \:  \: 
\tau \mapsto {\C}(1, \tau H, \tau ^{2}n). 
\end{equation} 
By sending   
$\infty \mapsto {\C}(0, 0, 1)$, 
the isomorphism $(\ref{eqn:tube domain realization})$ 
extends to the isomorphism 
\begin{equation}\label{eqn:extension of tube domain realization}
 {\mathbb H}^{\ast} :=  
 {\mathbb H} \cup {\Q} \cup \{ \infty \} 
\: \: \: \stackrel{\simeq}{\longrightarrow} \: \: \: 
\left(\Omega _{\widetilde{NS}(S)} ^{+} \right) ^{\ast} 
:=  \Omega _{\widetilde{NS}(S)}^{+} \cup 
\mathop{\bigcup}_{l \in \widetilde{NS}(S)\atop{(l, l)=0}} {\C}l .
\end{equation} 
Then 
the following diagram commutes:  
\begin{equation*}
\begin{array}{ccc}
  \Gamma _{0}(n)^{+}              &    \curvearrowright            &   {\mathbb H}^{\ast}    \\
     \downarrow                            &                                                 &    \downarrow                   \\
  \Gamma (S)^{+}                       &    \curvearrowright          &   \left( \Omega _{\widetilde{NS}(S)}^{+} \right) ^{\ast} .
\end{array}
\end{equation*}

In this way 
we obtain the desired isomorphism.

\begin{proposition}[cf.\cite{Do}]\label{Kahler moduli and Fricke curve}
The tube domain realization $(\ref{eqn:tube domain realization})$ 
induces the isomorphism 
\begin{equation}\label{eqn:Fricke modular curve and Kahler moduli}
\Gamma _{0}(n)^{+} \backslash {\mathbb H} 
\stackrel{\simeq}{\longrightarrow} 
\Gamma (S)^{+} \backslash \Omega _{\widetilde{NS}(S)}^{+}, 
\end{equation}
which extends naturally to 
the isomorphism between 
the compactifications of both sides.
\end{proposition}

The curve 
$\Gamma _{0}(n)^{+} \backslash {\mathbb H}$ 
is called 
the {\it Fricke modular curve\/}.

\subsection{The cusps of the Fricke modular curve}\label{sec5.2}

As is well-known, 
the set of the 
$\Gamma _{0}(n)$-cusps 
is naturally identified with the following finite set: 
\begin{equation}\label{eqn:cusps of Hecke modular curve}
\Bigl\{ 
\Bigl. \: (k, e) \: \Bigr| \: 
e \in {\Z}_{>0},  \: \: 
e | n ,  \: \: 
k \in ({\Z}/d{\Z})^{\times} \: {\rm where\/} \: d=(e, \frac{n}{e}) \: 
\Bigr\} . 
\end{equation}
The correspondence is as follows. 
For a rational number 
$\frac{s}{r} \in {\Q}$ 
with 
$r$ and $s$ coprime, 
we put  
$e:=(r, n)$ 
and 
$k:=\alpha s \in ({\Z}/d{\Z})^{\times}$ 
where 
$\alpha \in ({\Z}/n{\Z})^{\times}$ is such that 
$\alpha e = r \in {\Z}/n{\Z}$. 
Conversely, 
for a pair $(k, e)$ in the set $(\ref{eqn:cusps of Hecke modular curve})$  
we associate a rational number 
$\frac{\tilde{k}}{e} \in {\Q}$, 
where $\tilde{k} \in {\Z}$ is 
such that $\tilde{k} \equiv k \in {\Z}/d{\Z}$ 
and $(\tilde{k}, e)=1$.

The Fricke involution 
$(
\begin{smallmatrix}
0            &  -\sqrt{n}^{-1} \\
\sqrt{n}     &   0 
\end{smallmatrix}
)$
acts on the set $(\ref{eqn:cusps of Hecke modular curve})$
by 
$(k, e) \mapsto (-k, \frac{n}{e})$.   
Hence the set of the 
$\Gamma _{0}(n)^{+}$-cusps 
is identified with the following finite set: 
\begin{eqnarray*}\label{eqn:cusps of Fricke curve pre}
&   & 
         \Bigl\{ \Bigl. \:
                 (k, e) \: \: \Bigr| \: 
                                e \in {\Z}_{>0}, \: \: e|n, \: \: e>\frac{n}{e}, \: \: 
                                k \in ({\Z}/d{\Z})^{\times} \: {\rm where\/} \: d=(e, \frac{n}{e}) \: 
          \Bigr\} \\
& \bigcup &  
         \Bigl\{ \Bigl. \: 
                 (k', e) \: \: \Bigr| \: 
                                e \in {\Z}_{>0}, \: e^{2}=n, \:  
                                k' \in ({\Z}/e{\Z})^{\times}/ \{  \pm {\rm id\/} \: \} \:  
          \Bigr\} . 
\end{eqnarray*}
Writing $e=dr$ and $\frac{n}{e} = ds$, 
we obtain the following description.

\begin{proposition}\label{prop:cusps of Fricke modular curve} 
The set of the 
$\Gamma _{0}(n)^{+}$-cusps 
is identified with the following set:  
\begin{eqnarray}\label{eqn:cusps of Fricke curve final}
\mathcal{C}(n)^{+}    & :=  &   \bigcup_{{d^{2}|n}\atop{d^{2}<n}} 
                \Bigl\{ \Bigl. \:
                         (k, dr) \: \: \Bigr| \: 
                                       k \in ({\Z}/d{\Z})^{\times}, \: \: 
                                       \frac{n}{d^{2}}=rs, \: (r, s)=1,\:  r>s \:  
                \Bigr\} \nonumber\\
                 & \cup &   
                   \Bigl\{ \Bigl. \: 
                              (k', d) \: \: \Bigr| \: 
                                     d^{2}=n, \:  k' \in ({\Z}/d{\Z})^{\times}/ \{  \pm {\rm id\/} \: \} \:  
                   \Bigr\} . 
\end{eqnarray} 
Here $d$ runs over ${\Z}_{>0}$. 
\end{proposition}


On the other hand, 
the set ${\rm FM\/}^{d}(S)$ is described in \cite{Ma2},    
which we recall now briefly.  
Let $d$ be a positive integer satisfying $d^{2}|n$. 
For each $k \in ({\Z}/d{\Z})^{\times}$ 
choose a natural number $\tilde{k}$ 
such that $\tilde{k} \equiv k \in {\Z}/d{\Z}$ 
and $(\tilde{k}, 2n)=1$. 
When $d^{2}<n$, 
each element $\sigma \in ({\Z}/2{\Z})^{\tau (d^{-2}n)}$ 
corresponds to 
a pair $(r_{\sigma},  s_{\sigma})$ of coprime natural numbers  
satisfying 
$r_{\sigma}s_{\sigma} = d^{-2}n$. 
Let 
$\Sigma _{d} \subset ({\Z}/2{\Z})^{\tau (d^{-2}n)}$ 
be the subset 
consisting of those pairs $(r_{\sigma},  s_{\sigma})$ 
satisfying $r_{\sigma}>s_{\sigma}$. 
When $d^{2}=n$, 
consider the element $\sigma $ 
corresponding to the pair $(r_{\sigma}, s_{\sigma}) = (1, 1)$ 
and put $\Sigma _{d} := \{ \sigma \} $.

For a pair $(\sigma , k) \in \Sigma _{d} \times ({\Z}/d{\Z})^{\times}$  
we define 
\begin{equation*}\label{Mukai vectors}
v_{\sigma , k} := 
(dr_{\sigma}, \tilde{k}H, \tilde{k}^{2}ds_{\sigma}) \: 
\in I^{d}(\widetilde{NS}(S)). 
\end{equation*} 
Let 
$(M_{\sigma , k}, \alpha _{\sigma , k})$ 
be the moduli space of stable sheaves on $S$ 
with Mukai vector $v_{\sigma , k}$, 
where 
$\alpha _{\sigma , k}$ 
is the obstruction to the existence of a universal sheaf (\cite{Ca2}). 
We have the isomorphism 
\begin{equation*}
(M_{\sigma , k}, \alpha _{\sigma , k}) 
\simeq 
(S_{v_{\sigma , k}}, \alpha _{v_{\sigma , k}}), 
\end{equation*} 
where 
$(S_{v_{\sigma , k}}, \alpha _{v_{\sigma , k}})$ 
is the twisted $K3$ surface 
constructed from the cusp $[v_{\sigma , k}]$. 
By calculating the Hodge structures of $M_{\sigma , k}$ 
and applying the formula for $\# {\rm FM\/}^{d}(S)$, 
the following result is obtained in \cite{Ma2}.

\begin{proposition}[\cite{Ma2}]\label{Ma}
We have ${\rm FM\/}^{d}(S) = \emptyset $ unless $d^{2}|n$. 
Let $d^{2}|n$. 

$(1)$ If $d^{2}<n$, then 
\begin{equation*}
{\rm FM\/}^{d}(S) = 
            \Bigl\{ \Bigl. 
                            (M_{\sigma , k}, \alpha _{\sigma , k}) \: 
                    \Bigr|  \: 
                          (\sigma , k) \in  \Sigma _{d} \times ({\Z}/d{\Z})^{\times}
            \Bigr\} . 
\end{equation*}

$(2)$ Let $d^{2}=n$. 
Choose a set $\{ j \} \subset ({\Z}/d{\Z})^{\times}$ 
of representatives of the quotient set 
$({\Z}/d{\Z})^{\times} / \{ \pm {\rm id \/} \} $. 
Then 
\begin{equation*}
{\rm FM\/}^{d}(S) = 
            \Bigl\{ \Bigl. 
                            (M_{\sigma , k'}, \alpha _{\sigma , k'}) 
                    \: \Bigr| \: 
                          (\sigma , k') \in  \Sigma _{d} \times \{ j \}
            \Bigr\} . 
\end{equation*}
\end{proposition}

Via the tube domain realization  
$(\ref{eqn:extension of tube domain realization})$, 
we obtain the correspondence  
\begin{equation}\label{moduli and cusp}
(M_{\sigma , k}, \alpha _{\sigma , k})  \: \: \: 
\longleftrightarrow     \: \: \:
v_{\sigma , k}  \: \: \: 
             \stackrel{(\ref{eqn:extension of tube domain realization})}{\longleftrightarrow}     \: \: \: 
                                                              \frac{\tilde{k}}{dr_{\sigma}}    \: \: \:     
            \mapsto     \: \: \: 
                                 (k, dr_{\sigma}) \in \mathcal{C}(n)^{+}. 
\end{equation}              
Comparing 
Proposition \ref{prop:cusps of Fricke modular curve}
and 
Proposition \ref{Ma}, 
we now observe that 
$\mathop{\bigcup}_{d} {\rm FM\/}^{d}(S)$  
and  
$\mathcal{C}(n)^{+}$
corresponds bijectively, 
as expected. 
For an arbitrary rational number $\frac{a}{b} \in {\Q}$ 
with $a$ coprime to $b$,  
we can calculate 
the isomorphism class of the corresponding twisted FM partner 
with the aid of 
$(\ref{eqn:cusps of Fricke curve final})$ and $(\ref{moduli and cusp})$. 
In particular, 
we see that the $K3$ surface underlying the partner  
is determined only by the value $[b] \in {\Z}/n{\Z}$.    
In terms of a primitive isotropic Mukai vector 
$(r, kH, s)$, 
the $K3$ surface underlying the twisted moduli space 
is determined by the value 
\begin{equation*}
[(r, k)^{-1} r] \: \in {\Z}/n{\Z}.
\end{equation*}

Conversely, 
we may find the Mukai vectors $v_{\sigma , k}$ 
by referring to the cusps of the Fricke modular curve. 
If we resort to 
Theorem \ref{most general}, 
we reprove that the twisted moduli spaces $(M_{\sigma , k}, \alpha _{\sigma , k})$ 
represent ${\rm FM\/}^{d}(S)$, 
without calculating the twisted FM number of $S$ 
nor the Hodge structures of the moduli spaces $M_{\sigma , k}$.

In these ways,  
the classification of the twisted FM partners is 
simplified and is strengthened by using the Fricke modular curve.

\begin{example}
Let 
$S \subset {\proj}^{5}$ 
be a generic $K3$ surface of degree $8$ 
whose Picard group is generated by the class $H := \mathcal{O}(1)|_{S}$.    
It is well-known that $S$ is a $(2, 2, 2)$ complete intersection. 
The K\"ahler moduli 
$\mathcal{K} (S) \simeq\Gamma _{0}(4)^{+} \backslash {\mathbb H}$ 
has two cusps 
$[\infty]$ and $[\frac{1}{2}]$. 
The cusp $[\infty]$  corresponds to $S$ itself, 
while the cusps $[\frac{1}{2}]$ corresponds to 
the twisted moduli space $(M, \alpha )$ associated with the Mukai vector $(2, H, 2)$. 
The underlying $K3$ surface $M$ is of degree $2$, 
and can be realized as a double covering of  ${\proj}^{2}$ 
ramified over a sextic. 
\end{example}

\end{document}